\documentclass [11pt,twoside,a4paper]{article}
\usepackage{amsfonts}
\usepackage{amsmath}
\usepackage{amstext}
\usepackage{amssymb}
\usepackage{mathrsfs}
\usepackage{amscd}
\usepackage{xypic}
\usepackage{bbm}
\usepackage{mathabx}
\usepackage{epsf}              
\usepackage{graphicx}          
\usepackage{fancybox}          
\usepackage{color}             
\usepackage{fancyhdr}
\usepackage[hang,footnotesize]{caption2}  

\def\mathcal{\mathscr}
\newfont{\aaa}{cmb10 at 19pt}
\newfont{\bbb}{cmb10 at 11pt}
\newtheorem{thm}{Theorem}[section]
\newtheorem{cor}[thm]{Corollary}
\newtheorem{lem}[thm]{Lemma}
\newtheorem{prop}[thm]{Proposition}

\newtheorem{rem}[thm]{Remark}
\newtheorem{exm}[thm]{Example}
\def\le{\leqslant}
\def\leq{\leqslant}
\def\ge{\geqslant}

\newcommand{\beq}{\begin{equation}}
\newcommand{\eeq}{\end{equation}}
\newcommand{\bey}{\begin{eqnarray}}
\newcommand{\eey}{\end{eqnarray}}
\newcommand{\beyy}{\begin{eqnarray*}}
\newcommand{\eeyy}{\end{eqnarray*}}

\makeatletter
\def\@evenhead{
\vbox{\hbox to \textwidth {}{\hspace{0mm}{\footnotesize
\thepage}}{\hspace{8cm} {\footnotesize {Mu-Fa CHEN, Yu-Hui ZHANG}}} \protect\vspace{1truemm}\relax \hrule depth0pt
height0.15truemm width\textwidth}}
\def\@evenfoot{}
\def\@oddhead{\vbox{\hbox to \textwidth
{{\hspace{0cm}{\footnotesize Unified representation of formulas for single birth processes}\hfill{\footnotesize
\thepage}}\hspace{0mm}}{} \protect\vspace{1truemm}\relax\hrule
depth0pt height0.15truemm width\textwidth}}
\def\@oddfoot{}
\makeatother

\def\le{\leqslant}
\def\ge{\geqslant}
\def\leq{\leqslant}

\newcommand{\rf}[2]{[\ref{#1}; #2]}

\begin{document}

\thispagestyle{empty} \thispagestyle{fancy} {
\fancyhead[lO,RE]{\footnotesize  Front. Math. China 2014, 9(4): 761--796\\
DOI 10.1007/s11464-013-0186-5\\[3mm]
\includegraphics[0,-50][0,0]{11.bmp}}
\fancyhead[RO,LE]{\scriptsize \bf 
} \fancyfoot[CE,CO]{}}
\renewcommand{\headrulewidth}{0pt}


\setcounter{page}{1}
\qquad\\[8mm]

\noindent{\aaa{Unified representation of formulas for single birth processes}}\\[1mm]

\noindent{\bbb Mu-Fa CHEN,\quad Yu-Hui ZHANG}\\[-1mm]

\noindent\footnotesize{School of Mathematical Sciences, Beijing Normal University, Laboratory of Mathematics and Complex Systems, Ministry of Education,
Beijing 100875, China}\\[6mm]

\vskip-2mm \noindent{\footnotesize$\copyright$ Higher Education
Press and Springer-Verlag Berlin Heidelberg 2013} \vskip 4mm

\normalsize\noindent{\bbb Abstract}\quad Based on a new explicit representation of the solution to the Poisson equation with respect to single birth processes,
the unified treatment for various criteria on classical problems (including
uniqueness, recurrence, ergodicity, exponential ergodicity, strong ergodicity,
 as well as extinction probability etc.) for the processes are presented.\vspace{0.3cm}

\footnotetext{Received March 1, 2014; accepted April 5,
2014\\
\hspace*{5.8mm}Corresponding author: Yu-Hui ZHANG, E-mail: zhangyh@bnu.edu.cn}

\noindent{\bbb Keywords}\quad Single birth process, Poisson equation, uniqueness, recurrence, ergodicity, moments of return time\\
{\bbb MSC}\quad 60J60\\[0.4cm]

\setcounter{equation}{0}
\setcounter{section}{1}
\setcounter{thm}{0}

\noindent{\bbb 1\quad Introduction}\\[0.1cm]

\noindent
Consider a continuous-time homogeneous Markov chains $\{X(t):t\ge 0\}$, on a probability space $(\Omega, \mathscr F, \mathbb P)$,
with transition
probability matrix $P(t)=(p_{ij}(t))$ on a countable state space $\mathbb Z_+=\{0,1,2,\dots\}$.
We call $\{X(t):t\ge 0\}$ a single birth process if its transition rate (density) matrix $Q=(q_{ij}:i,j\in \mathbb Z_+)$
is irreducible and satisfies that $q_{i,i+1}>0,\, q_{i,i+j}=0$ for all $i\in \mathbb Z_+$ and $j\ge 2$.
Such a matrix $Q=(q_{ij})$ with $\sum_jq_{ij}=0$ for every $i$ (conservativity) is called a single birth $Q$-matrix. Refer to \cite{YC}.
In the literature, the single birth process is also called upwardly skip-free  process, or birth and death
process with catastrophes (cf. \cite{A91, B86, BGR} for instance).

The single birth process, as a natural extension of birth and death process which is a simplest $Q$-process (Markov chain),
has its own origins in practice, refer to the earlier papers \cite{B86, R, YC}, for instance.
The exit boundary of the process consists at most one single extremal point and so
the single birth process is nearly the largest class for which the explicit criteria on classical
problems can be expected. Actually, the study on the object is quite fruited and relatively completed (cf. \cite{C1, C2, C5, YC, ZJ, Z1}).
Based on this advantage, the single birth process becomes a fundamental comparison tool in studying
more complex processes, such as
infinite-dimensional reaction-diffusion processes. Refer to \rf{C1}{Chapters 3 and 4, Part I\!I\!I} and \cite{YC}.
Usually, the single birth process is non-symmetric and hence it is regarded as a representative one of the non-symmetric processes.
For non-symmetric processes, comparing with the symmetric ones, our knowledge is much limited, except for single birth processes
to which much results are known as just mentioned. Up to now, the known results are all presented in some recursive forms.
This paper introduces a single unified representation, as well as a unified treatment, of various formulas for single birth processes.

Throughout the paper, we consider only the single birth $Q$-matrix $Q=(q_{ij})$. Set
$q_i=-q_{ii}$ for each $i\in \mathbb Z_+$.
For a given function $c$ (to be fixed in this and the next sections, and then to be specified case by case), define an operator $\Omega$ as follows
$$\Omega g=Qg+cg, \quad \text{where}\quad (Qg)_i=\sum_{j}q_{ij}(g_j-g_i).$$
Clearly, if $c\le 0$, then $\Omega$ is an operator corresponding to a single birth
process with killing rates $(-c_i)$.

The following sequences are used throughout this paper.
\begin{gather}
\widetilde F_i^{(i)}=1,\qquad  \widetilde F_n^{(i)}=\frac 1{q_{n,n+1}}
\sum_{k=i}^{n-1}\tilde q_n^{(k)}\widetilde F_k^{(i)},\qquad n>i\ge 0,\label{ee0.01}\\
\tilde q_n^{(k)}=q_n^{(k)}-c_n:=\sum_{j=0}^kq_{nj}-c_n,\qquad 0\le k<n.\label{ee0.02}\end{gather}
Note that if $c\le 0$, then $\tilde q_n^{(k)}\ge 0$ and then
$\widetilde F_n^{(k)}\ge 0$ for every $n>k\ge 0$.
In what follows, we omit the superscript $\, \tilde{}\, $ everywhere in $\widetilde F$ and $\tilde q$ once $c_i\equiv 0$,
and often use the convention that $\sum_{\emptyset}=0$.

Here is the first of our main results.

\begin{thm}\label{t2.1}{ Given a single-birth $Q$-matrix $Q=(q_{ij})$ and functions $c$ and $f$,
the solution $g$ to the Poisson equation
\begin{equation}\label{ee0.1}
\Omega g=f
\end{equation}
has the following representation:
\begin{equation}\label{ee0.1-1}
g_n=g_0+\sum_{0\le k\le n-1}\sum_{0\le j\leq k}\frac{{\widetilde F}_k^{(j)}(f_j-c_jg_0)}{q_{j,j+1}},\qquad n\ge 0.
\end{equation}
In particular, the harmonic function $g$ of $\Omega$ $($i.e., $\Omega g=0)$ can be represented as
$$g_n=g_0\Bigg(1-\sum_{0\le k\le n-1}\sum_{0\le j\le k}\frac{\widetilde F_k^{(j)}c_j}{q_{j,j+1}}\Bigg),\qquad n\ge 0.$$

Conversely, for each boundary/initial value $g_0\in {\mathbb R}$, the function $(g_n)$ defined by $(\ref{ee0.1-1})$ is a solution
to $(\ref{ee0.1})$.
}\end{thm}

For single birth processes, almost all problems we concerned with are related to the solutions to some specific Poisson equation. Here, we unify these equations as (\ref{ee0.1}) with different functions $c$ and $f$ which are listed as follow.

\vskip0.3cm
\begin{center}{\renewcommand{\arraystretch}{1.6}
\begin{tabular}{|c|c|c|}
  \hline
  Problem                  & $c_i\in {\mathbb R}$                   & $f_i\in {\mathbb R}$
  \\ \hline
  Harmonic function               & $c_i\in {\mathbb R}$   &
  $f_i\equiv 0$ \\ \hline
  Uniqueness                & $c_i\equiv -\lambda<0$   & $f_i\equiv 0$ \\ \hline
  Recurrence                & $c_i\equiv 0$            & $f_i=q_{i0}(1-\delta_{i0})$\\ \hline
  Extinction/return probability    &$c_i\equiv 0$             & $f_i=q_{i0}(1-\delta_{i0})(g_0-1)$ \\ \hline
  Ergodicity                & $c_i\equiv 0$            & $f_i=q_{i0}(1-\delta_{i0})g_0-1$\\ \hline
  Strong ergodicity         & $c_i\equiv 0$            & $f_i=q_{i0}(1-\delta_{i0})g_0-1$\\ \hline
  Polynomial moment         & $c_i\equiv 0$            & $f_i^{(\ell)}$\\ \hline
  Exponential moment/ergodicity        & $c_i\equiv\lambda>0$     & $f_i=q_{i0}(1-\delta_{i0})(g_0-1)$ \\ \hline
  Laplace transform of return time & $c_i\equiv -\lambda<0$ & $f_i=q_{i0}(1-\delta_{i0})(g_0-1)$ \\ \hline
\end{tabular}}\end{center}
where $f_i^{(\ell)}=q_{ii_0}(1-\delta_{ii_0})g_{i_0}-\ell\mathbb E_i\sigma_{i_0}^{\ell-1}$.
\vskip0.32cm

We remark that in the two cases for ergodicity and strong ergodicity, even though the Poisson equation and the functions $c$ and $f$
are the same, but their solutions are required to be finite and bounded, respectively.

The paper is organized as follows. The proof of Theorem \ref{t2.1} is given in the next section,
 using a lemma on the representation of solution to a class of linear equations. Then,
Sections 3--7 are devoted, respectively, to the criteria on the problems listed in the table above, and related problems to be specific subsequently.
Roughly speaking, the unified treatment
presented in the paper consists of the following three steps.
\begin{itemize}\setlength{\itemsep}{-0.6ex}
\item[(a)] Find out the Poisson equation corresponding to the problem we are interested in.
\item[(b)] Apply Theorem \ref{t2.1} to get the solution to the Poisson equation.
\item[(c)] Work out a criterion for the problem using the solution obtained in (b).
\end{itemize}
Step (a) is more or less known from the previous study; step (b) is now automatic; hence, our main
work is spent on step (c).

For the reader's convenience, several key formulas used often in the proofs are collected into an Appendix in a single page which consists the last page of the paper (so that it can be printed out separately).
\\[4mm]

\setcounter{equation}{0}
\setcounter{section}{2}
\setcounter{thm}{0}
\noindent{\bbb 2\quad The Poisson equation}\\[0.1cm]

\noindent
In this section, we consider the solutions of the Poisson equation (\ref{ee0.1}) for single birth processes.
Let us begin with a simple result for the solution to a class of linear equations.

\begin{lem}\label{l00}{ For given real numbers $(\alpha_{nk})_{n-1\ge k\ge 0}$ and $(f_n)_{n\ge 0}$, the solution $(g_n)_{n\ge 0}$ to the recursive
inhomogeneous equations
\begin{equation}\label{e0.01}
g_n=\sum_{0\le k\le n-1}\alpha_{nk}g_k+f_n,\qquad n\ge 0
\end{equation}
can be represented as
\begin{equation}\label{e0.02}
g_n=\sum_{0\le k\le n}\gamma_{nk}f_k,\qquad n\ge 0,
\end{equation}
where for fixed $k\ge 0$, $(\gamma_{nk})_{n\ge k}$ with $\gamma_{kk}=1$ is the solution to the recursive equations
\begin{equation}\label{e0.03}\gamma_{nk}=\sum_{k\le j\le n-1}\alpha_{nj}\gamma_{jk},\qquad n>k.
\end{equation}
}\end{lem}

\noindent{\it Proof}\quad Use induction. For $n=0$, we have
$$g_0=f_0=\gamma_{00}f_0=\sum_{0\le k\le 0}\gamma_{0k}f_k.$$
Assume that (\ref{e0.02}) holds for all $n\le m$. When $n=m+1$, from (\ref{e0.01}), we see that
$$\aligned
g_{m+1}&=\sum_{0\le k\le m}\alpha_{m+1,k}\,g_k +f_{m+1}=
\sum_{0\le k\le m}\alpha_{m+1,k}\sum_{0\le\ell\le k}\gamma_{k\ell}\,f_\ell+f_{m+1}\\
&=\sum_{0\le\ell\le m}\Bigg(\sum_{\ell\le k\le m}\alpha_{m+1,k}\gamma_{k\ell}\Bigg)\,f_\ell
+f_{m+1}=\sum_{0\le\ell\le m}\gamma_{m+1,\ell}\,f_\ell+f_{m+1}\\
&=\sum_{0\le\ell\le m+1}\gamma_{m+1,\ell}\,f_\ell.
\endaligned$$
Hence, (\ref{e0.02}) holds for $n=m+1$. By induction, the representation (\ref{e0.02}) holds for all $n\ge 0$.
\hfill$\square$

\vskip0.2cm
Note that the coefficients $(\alpha_{nk})$ are often fixed and so are $(\gamma_{nk})$.
Then Lemma \ref{l00} says that once replacing $(\alpha_{nk})$ by
$(\gamma_{nk})$, the solution to (\ref{e0.01}) has a complete representation (\ref{e0.02}),
mainly in terms of the inhomogeneous term $(f_n)$ in (\ref{e0.01}).

Without condition $\gamma_{kk}=1$,  equation (\ref{e0.03}) is clearly homogeneous. However,
it becomes inhomogeneous under condition $\gamma_{kk}\ne 0$ (then one may assume that $\gamma_{kk}=1$):
$$\gamma_{nk}=\sum_{k+1\le j\le n-1}\alpha_{nj}\gamma_{jk}+\alpha_{nk}\gamma_{kk},\qquad n\ge k+1$$
provided $\alpha_{k+1, k}\ne 0$. Otherwise, once $\alpha_{k+1, k}=0$, by induction, we actually
have $\gamma_{nk}=0$ for all $n\ge k+1$. Thus, under condition $\gamma_{kk}=1$, by Lemma \ref{l00}
(for fixed $k$),
we have the following alternative representation of $(\gamma_{nk})$:
$$\gamma_{nk}=\sum_{k+1\le j\le n}\gamma_{nj}\alpha_{jk},\qquad n\ge k+1.$$

In what follows, we will use the following variant of Lemma \ref{l00}. Replacing the initial
$0$ by $i$ and the coefficient $(\alpha_{nk})$ by $(\alpha_{nk}\beta_k)$,
respectively, for some non-zero sequence
$(\beta_n)$, and set $h_n=g_n/\beta_n\,(n\ge i)$, we obtain the following result.

\begin{cor}\label{c0}{The solution $(h_n)_{n\ge i}$ to the recursive equations
\begin{equation}\label{e0.1}
h_n=\frac 1{\beta_n}\Bigg(\sum_{i\le k\le n-1}\alpha_{nk}h_k+f_n\Bigg),\qquad n\ge i
\end{equation}
can be represented as
\begin{equation}\label{e0.2}
h_n=\sum_{i\le k\le n}\frac{\gamma_{nk}}{\beta_k}f_k,\qquad n\ge i,
\end{equation}
where for each fixed $i$, $(\gamma_{ni})_{n\ge i}$ with $\gamma_{ii}=1$ is the solution to the equations
$$\gamma_{ni}=\frac 1{\beta_n}\sum_{i\le k\le n-1}\alpha_{nk}\gamma_{ki},\qquad n> i.$$
Equivalently,
\begin{equation}\label{e0.3}
\gamma_{ii}=1,\qquad\gamma_{ni}=\sum_{i+1\le k\le n}\frac{\gamma_{nk}}{\beta_k}\alpha_{ki},\qquad n\ge i+1.
\end{equation}
}\end{cor}

Specifying $\beta_n=q_{n,n+1}$ and $\alpha_{nk}={\tilde q}_n^{(k)}$ in Corollary \ref{c0}
and using the successive formula of ${\widetilde F}_n^{(k)}$ defined in (\ref{ee0.01}),
we obtain the following result.

\begin{cor}\label{c0-1}{ For given $f$, the sequence $(h_n)$ defined successively by
$$h_n=\frac{1}{q_{n, n+1}}\bigg(f_n+\sum_{i\le k \le n-1}{\tilde q}_n^{(k)} h_k\bigg),\qquad n\ge i$$
has an unified expression as follows
$$h_n=\sum_{k=i}^n \frac{{\widetilde F}_n^{(k)}}{q_{k, k+1}} f_k,\qquad n\ge i. $$
In particular, the sequence $\big({\widetilde F}_n^{(k)}\big)$ defined in $(\ref{ee0.01})$ has the following expression
\begin{equation}\label{e0.4-0}{\widetilde F}_i^{(i)}=1,\qquad  {\widetilde F}_n^{(i)}=\sum_{k=i+1}^n
\frac{{\widetilde F}_n^{(k)}{\tilde q}_k^{(i)}}{q_{k, k+1}},\quad n\ge i+1. \end{equation}
}\end{cor}

Before moving further, let us mention a comparison result for different $\gamma_{nj}$, which may be useful elsewhere but not
in this paper.

\begin{prop}\label{p1}{ For each triple $n\ge i> j$, the following assertion holds:
\begin{equation}\label{e0.4}
\gamma_{nj}=\sum_{i\le k\le n}\frac{\gamma_{nk}}{\beta_k}\sum_{j\le\ell\le i-1}\alpha_{k\ell}\gamma_{\ell j}.
\end{equation}
Furthermore, if $\alpha_{nk}\ge 0$ and $\beta_n>0$ for all $n>k$, then $\gamma_{ni}\gamma_{ij}\le \gamma_{nj}$ for all $n\ge i\ge j$.
}\end{prop}

\noindent{\it Proof\quad} The first assertion is simply a consequence of Corollary \ref{c0}. In fact, for fixed $i>j$, take
$$f_n=\sum_{j\le\ell\le i-1}\alpha_{n\ell}\,\gamma_{\ell j},\qquad n\ge i.$$
Then
$$\gamma_{nj}=\frac 1{\beta_n}\Bigg[\sum_{i\le\ell\le n-1}\!\alpha_{n,\ell}\,\gamma_{\ell j}+\!\!\sum_{j\le\ell\le i-1}\!
\alpha_{n\ell}\,\gamma_{\ell j}\Bigg]\!=\frac 1{\beta_n}\Bigg[\sum_{i\le\ell\le n-1}\!\alpha_{n\ell}\,\gamma_{\ell j}+f_n\Bigg],\quad n\ge i.
$$
Hence, by Corollary \ref{c0}, we get
$$\gamma_{nj}=\sum_{i\le k\le n}\frac{\gamma_{nk}}{\beta_k}f_k=\sum_{i\le k\le n}\frac{\gamma_{nk}}{\beta_k}\sum_{j\le\ell\le i-1}\alpha_{k\ell}
\,\gamma_{\ell j},\qquad n\ge i.$$

If $\alpha_{nk}\ge 0$ and $\beta_n>0$ for all $n$ and $k$, then from (\ref{e0.4}), it follows that for all $n> i>j$,
$$\gamma_{nj}=\gamma_{ni}\gamma_{ij}+\sum_{i+1\le k\le n}\frac{\gamma_{nk}}{\beta_k}\sum_{j\le \ell\le i-1}\alpha_{k\ell}\gamma_{\ell j}\ge \gamma_{ni}\gamma_{ij}.$$
In the cases of $n=i$ or $i=j$, the conclusion is trivial.\hfill$\square$

\vskip0.2cm
Now we turn to prove our first result.

\vskip0.2cm
\noindent{\it Proof of Theorem $\ref{t2.1}$}\quad For each $i\ge 0$, we have
\begin{align}
(\Omega g)_i&=q_{i,i+1}(g_{i+1}-g_i)-\sum_{0\le j\le i-1} q_{ij}\sum_{k=j}^{i-1}(g_{k+1}-g_k)+c_ig_i\nonumber\\
&=q_{i,i+1}(g_{i+1}-g_i)-\sum_{0\le k\le i-1}\sum_{j=0}^{k}q_{ij}(g_{k+1}-g_k)+c_ig_i\nonumber\\
&=q_{i,i+1}(g_{i+1}-g_i)-\sum_{0\le k\le i-1}\Bigg(\sum_{j=0}^{k}q_{ij}-c_i\Bigg)(g_{k+1}-g_k)+c_ig_0\nonumber\\
&=q_{i,i+1}(g_{i+1}-g_i)-\sum_{0\le k\le i-1}\tilde q_i^{(k)}(g_{k+1}-g_k)+c_ig_0.
\label{e0.4-1}\end{align}
Denote $g_{k+1}-g_k$ by $w_k$ for $k\ge 0$. Then
$$(\Omega g)_i=q_{i,i+1}w_i-\sum_{0\le k \le i-1}\tilde q_{i}^{(k)}w_k+c_ig_0,\qquad i\ge 0.$$
Now we rewrite the Poisson equation (\ref{ee0.1}) as
$$w_i=\frac 1{q_{i,i+1}}\Bigg(\sum_{0\le k \le i-1}\tilde q_{i}^{(k)}w_k+\tilde f_i\Bigg),\qquad i\ge 0,$$
where $\tilde f_i=f_i-c_ig_0$ for $i\ge 0$. By Corollary \ref{c0-1}, we obtain
$$w_i=\sum_{j=0}^i\frac{\widetilde F_i^{(j)}\tilde f_j}{q_{j,j+1}},\qquad i\ge 0.$$
So the solution of the Poisson equation (\ref{ee0.1}) satisfies
$$g_i=g_0+\sum_{k=0}^{i-1}w_k=g_0+\sum_{k=0}^{i-1}\sum_{j=0}^k\frac{\widetilde F_k^{(j)}\tilde f_j}{q_{j,j+1}},\qquad i\ge 1.$$
The first assertion is proven. The second assertion is simply a consequence of the first one.

To prove the last assertion of the theorem, noting that by (\ref{ee0.1-1}), we have
$$g_{n+1}-g_n=\sum_{j=0}^n\frac{{\widetilde F}_n^{(j)}(f_j-c_jg_0)}{q_{j,j+1}},\qquad n\ge 0.$$
Thus, from (\ref{e0.4-1}), it follows for each $i\ge 0$ that
$$\aligned
(\Omega g)_i&= q_{i, i+1} \sum_{j=0}^i\frac{{\widetilde F}_i^{(j)}(f_j-c_jg_0)}{q_{j,j+1}}
-\sum_{0\le k\le i-1}{\tilde q}_i^{(k)}  \sum_{j=0}^k\frac{{\widetilde F}_k^{(j)}(f_j-c_jg_0)}{q_{j,j+1}}+c_i g_0.
\endaligned$$
Because (by exchanging the order of sums and using (\ref{ee0.01}))
$$\aligned
\sum_{0\le k\le i-1}{\tilde q}_i^{(k)}  \sum_{j=0}^k\frac{{\widetilde F}_k^{(j)}(f_j-c_jg_0)}{q_{j,j+1}}
&=\sum_{0\le j\le i-1}\frac{f_j-c_jg_0}{q_{j,j+1}}\sum_{k=j}^{i-1} {\tilde q}_i^{(k)} {\widetilde F}_k^{(j)}\\
&=q_{i, i+1}\sum_{0\le j\le i-1}\frac{{\widetilde F}_i^{(j)}(f_j-c_jg_0)}{q_{j,j+1}},
\endaligned$$
we obtain $\Omega g =f$ as required.
\hfill$\square$

\begin{rem}{\rm  $(1)$ One may obtain $\big({\tilde q}_n^{(k)}, {\widetilde F}_n^{(k)}\big)$
from $\big({q}_n^{(k)}, {F}_n^{(k)}\big)$ easily replacing the
original $Q=(q_{ij})$ by ${\widetilde Q}=({\tilde q}_{ij})$:
$$\begin{cases}
{\tilde q}_{i0}= q_{i0}-c_i,\\
{\tilde q}_{ij}=q_{ij}, \quad j\ne 0,\; i\in E.
\end{cases}$$
In other words, only the first column of $Q=(q_{ij})$ is modified.
Then the original Poisson equation $\Omega g=f$ can be rewritten as
${\widetilde Q} g=\tilde f$ with ${\tilde f}_i= f_i-c_i g_0$.

$(2)$ Alternatively, one may enlarge the space $E$ by adding a point, say $-1$
for instance. Then introduce suitable ${\bar q}_{-1, i}$, ${\bar q}_{i, -1}$, ${\bar g}_{-1}$,
and ${\bar f}_{-1}$, so that ${\widebar Q}|_E=Q$, ${\bar g}|_E=g$, and ${\bar f}|_E=f$.
In this way, one may rewrite $\Omega g=f$ on $E$ as
${\widebar Q}{\bar g}= {\bar f}$ on $E\cup \{-1\}$.

$(3)$ To solve the Poisson equation, in view of $(\ref{e0.4-1})$, even for the simplest birth--death type, once $c$ appears,
it is necessary to go out to the larger class of single birth one, one can not just stay within the class of birth--death processes.
Actually, this observation is crucial to solve the Open Problem 9.13 in \cite{C6}. Refer to \rf{CZX}{Theorem 2.6}.
}\end{rem}

For the remainder of this section, we consider only the processes on a finite state space
$\{0, 1,\dots, N\}$. Note that here the rate $q_{N, N+1}$ is
not defined (or setting to be zero), but we allow $c_N\ne 0$. Hence ${\widetilde F}_n^{(k)}$ is
defined up to $n=N-1$ only. The next result is a localized version of Theorem \ref{t2.1}

\begin{prop}\label{p2.1}\!{  Given a single-birth $Q$-matrix
$(q_{ij})$ and a function $c$ on the finite state space $\{0, 1,\dots, N\}$ $(N\ge 1)$,
the following assertions hold.

$(i)$ The solution of the Poisson equation $\Omega g=f$ has the following form:
\begin{equation}\label{e2.2}
g_n=g_0
+\sum_{0\le k\le n-1}\sum_{0\le j\le k}\frac{{\widetilde F}_k^{(j)}(f_j-c_j g_0)}{q_{j,j+1}},\qquad 0\le n\le N,
\end{equation}
with boundary condition
$$c_Ng_0=\sum_{k=0}^{N-1}{\tilde q}_N^{(k)}\sum_{j=0}^k \frac{{\widetilde F}_k^{(j)}(f_j-c_j g_0)}{q_{j,j+1}}+f_N. $$

$(ii)$ Let $c\le 0$. Then the harmonic equation $\Omega g=0$ has only the trivial solution $g_i\equiv 0$
 iff there exists some $c_i< 0$.

$(iii)$ The unique solution $g$ to the equation $\Omega g|_{\{0,1,\dots,N-1\}}=0$ $($locally harmonic$)$ with $g_0=1$
is as follows:
\begin{equation}\label{e2.6}
g_n=1-\sum_{0\le k\le n-1}\sum_{0\le j\le k}\frac{{\widetilde F}_k^{(j)}c_j}{q_{j,j+1}},\qquad 0\le n\le N
\end{equation}
which is increasing once $c\le 0$.
}\end{prop}

\noindent{\it Proof}\quad (a) The proof is nearly the same as the one of Theorem \ref{t2.1}, except
we have to take care for the boundary at $N$. By (\ref{e0.4-1}), for $0\le i\le N-1$, we have
$$(\Omega g)_i=q_{i,i+1}(g_{i+1}-g_i)-\sum_{0\le k\le i-1}{\tilde q}_{i}^{(k)}(g_{k+1}-g_k)+c_ig_0.$$
Denote $g_{k+1}-g_k$ by $w_k$ for all $0\le k<N$.
Then
$$\aligned &(\Omega g)_i=q_{i,i+1}w_i-\sum_{0\le k\le i-1}{\tilde q}_{i}^{(k)}w_k+c_ig_0,\qquad 0\le i<N;\\
& (\Omega g)_N=-\sum_{k=0}^{N-1}{\tilde q}_N^{(k)}w_k+c_Ng_0.\endaligned$$
Rewrite the Poisson equation as
\begin{equation}\label{e2.3}
w_i=\frac 1{q_{i,i+1}}\Bigg(\tilde f_i+\sum_{0\le k\le i-1}\tilde q_{i}^{(k)}w_k\Bigg),\qquad 0\le i< N,
\end{equation}
where $\tilde f_i=f_i-c_ig_0$ for all $0\le i\le N$.
By Corollary \ref{c0-1}, we get
\begin{equation}\label{e2.4}
w_i=\sum_{j=0}^i\frac{{\widetilde F}_i^{(j)}\tilde f_j}{q_{j,j+1}},\qquad 0\le i< N.
\end{equation}
So the solution of the Poisson equation satisfies
$$g_i=g_0+\sum_{k=0}^{i-1}w_k=g_0+\sum_{k=0}^{i-1}\sum_{j=0}^k\frac{{\widetilde F}_k^{(j)}\tilde f_j}{q_{j,j+1}},\qquad 1\le i\le N.$$
Combining this with the boundary condition $(\Omega g)_N=f_N$ and (\ref{e2.4}), we obtain the first assertion.

(b) We have just seen that the harmonic solution $g$ satisfies
\begin{equation}\label{e2.5}
g_n=g_0\Bigg(1-\sum_{k=0}^{n-1}\sum_{j=0}^k\frac{{\widetilde F}_k^{(j)}c_j}{q_{j,j+1}}\Bigg),\qquad 1\le n\le N.
\end{equation}
and
$$g_0\Bigg(c_N+\sum_{k=0}^{N-1}\tilde q_N^{(k)}\sum_{j=0}^k\frac{{\widetilde F}_k^{(j)}c_j}{q_{j,j+1}}\Bigg)=0.$$
When $c\le 0$, by irreducibility, we have not only ${\tilde q}_N^{(N-1)}>0$ but also ${\widetilde F}_{N-1}^{(j)}>0$ for every $j: 0\le j\le N-1$.
Hence, if there exists some $c_i<0$, then we must have $g_0=0$ by the last equation.
Furthermore, by (\ref{e2.5}), we indeed have $g\equiv 0$.

Conversely, if $c_i\equiv 0$, then every constant function $g\ne 0$ is a solution to the equation $\Omega g=0$. Hence the harmonic
function $g$ can be non-trivial.

(c) To prove the third assertion, based on the second one, we have to use a smaller space $\{0,1,\ldots, N-1\}$ instead of the
original $\{0,1,\ldots, N\}$ to avoid the trivial solution. The assertion now follows from (\ref{e2.5}).\hfill$\square$

\vskip0.2cm

The next result is exceptional of the paper. Instead of single birth, we consider single death processes on a finite state space.
The result may be regarded as a dual of Proposition \ref{p2.1}. It indicates that a large parts of the study in the paper is
meaningful for the single death processes, but we will not go
to the details here.

A matrix $Q=(q_{ij})$ is called of single death if $q_{i,i-j}>0$
iff $j=1$ for $i\ge 1$.

\begin{prop}\label{p2.2}{  Given a single death $Q$-matrix $Q=(q_{ij})$ and a function $(c_i)$
on the finite state space $\{0, 1,\dots, N\}$, define
${\tilde q}_n^{(k)}=\sum_{j=k}^N q_{nj}-c_n$ for $k> n$ and
$${\widetilde F}_i^{(i)}=1,\qquad  {\widetilde F}_n^{(i)}=\frac 1{q_{n,n-1}}
\sum_{k=n+1}^{i}{\tilde q}_n^{(k)} {\widetilde F}_k^{(i)},\qquad 1\le n<i.$$
Then

$(i)$ the solution $g$ to the Poisson equation $\Omega g=f$ has the following representation:
$$g_n=g_N
+\sum_{n+1\le k\le N}\sum_{k\le j\le N}\frac{{\widetilde F}_k^{(j)}(f_j-c_j g_N)}{q_{j,j-1}},\qquad 0\le n\le N$$
with boundary condition
$$c_0^{} g_N^{}= \sum_{k=1}^N {\tilde q}_0^{(k)}\sum_{j=k}^N \frac{{\widetilde F}_k^{(j)}(f_j-c_j g_N)}{q_{j,j-1}}+f_0.$$
$(ii)$ The unique solution with $g_N^{}=1$ to equation $Qg|_{\{1,2,\dots,N\}}=0$ is as follows:
$$
g_n=1-\sum_{n+1\le k\le N}\sum_{k\le j\le N}\frac{{\widetilde F}_k^{(j)}c_j}{q_{j,j-1}}\qquad (0\le n\le N)
$$
which is decreasing in $n$ once $c\le 0$.
}\end{prop}

\noindent{\it Proof}\quad For $1\le i\le N$, we have
$$\aligned
(\Omega g)_i&=q_{i,i-1}(g_{i-1}-g_i)+\sum_{i+1\le j\le N}q_{ij}\sum_{k=i+1}^{j}(g_{k}-g_{k-1})+c_ig_i\\
&=q_{i,i-1}(g_{i-1}-g_i)+\sum_{i+1\le k\le N}\sum_{j=k}^{N}q_{ij}(g_{k}-g_{k-1})+c_ig_i\\
&=q_{i,i-1}(g_{i-1}-g_i)-\sum_{i+1\le k\le N}{\tilde q}_{i}^{(k)}(g_{k-1}-g_{k})+c_ig_N.
\endaligned$$
Denote $g_{k-1}-g_k$ by $w_k$ for all $1\le k\le N$. Then
\begin{gather}
\quad (\Omega g)_i=q_{i,i-1}w_i-\sum_{i+1\le j\le N}{\tilde q}_{i}^{(k)}w_k+c_ig_N,\qquad 1\le i\le N;\nonumber\\
 (\Omega g)_0=-\sum_{k=1}^{N}{\tilde q}_{0}^{(k)}w_k+c_0g_N.\nonumber\end{gather}
Now we rewrite the Poisson equation as
$$w_i=\frac 1{q_{i,i-1}}\Bigg(\tilde f_i+\sum_{i+1\le j\le N}q_{i}^{(k)}w_k\Bigg),\qquad 1\le i\le N,$$
where $\tilde f_i=f_i-c_ig_N$ for all $0\le i\le N$.
As an analogue of Corollary \ref{c0-1}, by induction, we can verify that
$$w_i=\sum_{j=i}^N\frac{\widetilde F_i^{(j)}\tilde f_j}{q_{j,j-1}},\qquad 1\le i\le N.$$
From the argument above, it follows immediately that
$$g_i=g_N+\sum_{k=i+1}^{N}w_k=g_N+\sum_{i+1\le j\le N}\sum_{k\le j\le N}\frac{ \widetilde F_k^{(j)}\tilde f_j}{q_{j,j-1}},\qquad  0\le i\le N-1.$$
Combining this with the boundary condition $(\Omega g)_0=f_0$, we finish the proof of the first assertion.
The second assertion is derived from the first one immediately. \hfill$\square$
\\[4mm]

\setcounter{equation}{0}
\setcounter{section}{3}
\setcounter{thm}{0}
\noindent{\bbb 3\quad Uniqueness}\\[0.1cm]

\noindent
Starting from this section, we handle with the problems
for single birth processes,  listed at the beginning of the paper. First, we study the uniqueness problem.
To do so, we need a sequence $({\widetilde m}_n)$(to be used often subsequently) :
\begin{equation}\label{e3.4-00}{\widetilde m}_0=\frac{1}{q_{01}},\quad
{\widetilde m}_n=\frac 1{q_{n,n+1}}
\bigg(1+\sum_{k=0}^{n-1}{\tilde q}_n^{(k)}{\widetilde m}_k\bigg),\qquad n\ge 1.\end{equation}
By Corollary \ref{c0-1}, we have
\begin{equation}\label{e3.4-0}{\widetilde m}_n
=\sum_{k=0}^n\frac{{\widetilde F}_n^{(k)}}{q_{k,k+1}},\qquad n\ge 0.\end{equation}
Again, we omit the superscript $\,\tilde{}\,$ everywhere in $\widetilde m$,
$\widetilde F$, and $\tilde q$ once $c_i\equiv 0$.
The following criterion is taken from \cite{C1, YC, ZJ}.

\begin{prop}\label{pr3-1}{  Corresponding to a given single birth $Q$-matrix $Q=(q_{ij})$ $($conservative$)$, the process is unique $($non-explosive$)$ iff
$\sum_{n=0}^\infty m_n=\infty$.
}\end{prop}

\noindent{\it Proof}\quad By \rf{C1}{Theorems 2.47 and 2.40}, the single birth process is unique iff the solution $(u_i)$ to the equation
\begin{equation}\label{e3.4}(\lambda +q_i)u_i=\sum_{j\ne i}q_{ij}u_j, \qquad i\ge 0;\qquad u_0=1
\end{equation}
is unbounded for some (equivalently for all) $\lambda>0$.
Rewrite (\ref{e3.4}) as
$$\Omega u=Qu- \lambda u=0;\qquad u_0=1.$$
Applying Theorem \ref{t2.1} to $c_i\equiv -\lambda$ and $f_i\equiv 0$, we obtain the unique solution:
$$u_n=1+\lambda\sum_{0\le k \le n-1}\sum_{j=0}^k\frac{\widetilde F_k^{(j)}}{q_{j,j+1}}=1+\lambda\sum_{0\le k \le n-1}\widetilde m_k,\qquad n\ge 0.$$
Clearly, $u_n$ is increasing in $n$ and then is unbounded iff $\sum_n {\widetilde  m}_n=\infty$. Thus,
it remains to show that $\sum_n {\widetilde  m}_n=\infty$
iff $\sum_n {m}_n=\infty$. Combining ${\widetilde  m}_n$ with $m_n$, it is clear that
$${\widetilde  m}_n =\sum_{j=0}^n\frac{\widetilde F_n^{(j)}}{q_{j,j+1}}\;\;\big\downarrow \;\;
\sum_{k=0}^n\frac{F_n^{(k)}}{q_{k,k+1}}=m_n\qquad
\text{as $\lambda\downarrow 0$},$$
since
$${\tilde q}_n^{(k)}={q}_n^{(k)}+\lambda \;\;\downarrow \;\; {q}_n^{(k)} \qquad
\text{as $\lambda\downarrow 0$}.$$
This already shows that the condition $\sum_n m_n=\infty$ is sufficient. It is nearly necessary
since the conclusion does not depend on $\lambda>0$, except there is a jump from $\lambda>0$
to $\lambda=0$. Hopefully, we
have thus seen some advantage of Theorem \ref{t2.1}, even though there is
still a distance to prove the necessity.

Actually, there are several ways to prove the equivalence
$$\sum_n {\widetilde  m}_n=\infty\;\text{for a fixed }\lambda>0 \Longleftrightarrow \sum_n {m}_n=\infty.$$
From now on, for simplicity, assume that $\lambda=1$.

(a) Observing that corresponding to the sequence $(\widetilde m_n)$, the operator is $\Omega=Q-I$
which may be regarded as a bounded perturbation of the original operator $Q$. Since these two operators
are zero-exit or not simultaneously, the equivalence above holds.

(b) In the original proof (cf. \rf{C1}{Proof of Theorem 3.16}), it was proved that $u_n$ is unbounded iff
$\sum_n m_n=\infty$. Combining this with what proved above, we obtain the required equivalence.

(c) Here is a more direct proof. The idea comes from \cite{Z4}.

Assume that $\sum_{k=0}^\infty\widetilde m_k=\infty$. If $\sum_{k=0}^\infty m_k<\infty$,
then there exists $N_0$ large enough such that
for all $n\ge N_0$,
$${\widetilde M}_n:=\sum_{k=0}^n\widetilde m_k>1\quad\text{and}\quad K:=2\sum_{k=N_0+1}^\infty m_k<1.$$
We now prove that for each $n>N_0$,
\begin{equation}\label{e3.5}
\widetilde m_k\le 2 m_k {\widetilde M}_{n-1},\qquad 0\le k\le n.
\end{equation}
Since $\widetilde m_0=m_0$ and ${\widetilde M}_{n-1}>1$ (due to the fact that $n-1\ge N_0$), (\ref{e3.5})
holds in the case of $k=0$. Assume that (\ref{e3.5}) holds up to $k=\ell-1<n$. Then,
$$\aligned
\widetilde m_\ell&=\frac 1{q_{\ell,\ell+1}}\Bigg(1+\sum_{k=0}^{\ell-1}q_\ell^{(k)}\widetilde m_k+\sum_{k=0}^{\ell-1}\widetilde m_k\Bigg)\quad\text{\rm (since $\lambda=1$)}\\
&\le \frac 1{q_{\ell,\ell+1}}\Bigg(1+\sum_{k=0}^{\ell-1}q_\ell^{(k)}2m_k{\widetilde M}_{n-1}+{\widetilde M}_{\ell-1}\Bigg)\quad\text{\rm (by assumption)}\\
&\le \frac 1{q_{\ell,\ell+1}}\Bigg(1+\sum_{k=0}^{\ell-1}q_\ell^{(k)}m_k\Bigg)2{\widetilde M}_{n-1}\\
&=2 m_\ell {\widetilde M}_{n-1}.
\endaligned$$
So (\ref{e3.5}) holds when $k=\ell$. By induction, we know that (\ref{e3.5}) holds for every $k: 0\le k\le n$.
Now, for each $n>N_0$, we have
$${\widetilde M}_n={\widetilde M}_{N_0}+\sum_{k=N_0+1}^n\widetilde m_k\le {\widetilde M}_{N_0}
+\sum_{k=N_0+1}^n 2 m_k{\widetilde M}_{n-1}\le {\widetilde M}_{N_0}+K{\widetilde M}_{n-1}.$$
Furthermore, we have
$$\aligned
{\widetilde M}_n &\le {\widetilde M}_{N_0}\big(1+K+\cdots+K^{n-N_0-1}\big)+K^{n-N_0}{\widetilde M}_{N_0}\\
&=\frac{{\widetilde M}_{N_0}(1-K^{n-N_0})}{1-K}+K^{n-N_0}{\widetilde M}_{N_0}.
\endaligned$$
Thus, as $n\to\infty$, we would have $\infty\le {{\widetilde M}_{N_0}}/{(1-K)}$  which is a contradiction. Hence,
once $\sum_{k=0}^\infty {\widetilde m}_k=\infty$,
we should also have $\sum_{k=0}^\infty m_k=\infty$.

We have therefore completed the proof of the equivalence mentioned above.

\hfill$\square$

\vskip0.2cm
To conclude this section, we mention that the uniqueness problem for the single birth $Q$-matrix
with absorbing set $H=\{0,1,\ldots, N\}\,(N<\infty)$
can be dealt with by the same approach. Refer to \rf{C1}{Theorem 3.16} and \cite{WZ}.
\\[4mm]

\setcounter{equation}{0}
\setcounter{section}{4}
\setcounter{thm}{0}
\noindent{\bbb 4\quad Recurrence and extinction/return probability}\\[0.1cm]

\noindent
For the recurrence, the following criterion is taken from \rf{C1}{Theorem 4.52\,(1)} and \cite{YC}.

\begin{prop}\label{pr}{  Assume the single birth $Q$-matrix $Q=(q_{ij})$ is non-explosive and irreducible. Then the process is recurrent iff $\sum_{n=0}^\infty F_n^{(0)}=\infty$, where $\big(F_n^{(i)}\big)$ was defined in $(\ref{ee0.01})$ by setting $c_i\equiv 0$.
}\end{prop}

\noindent{\it Proof}\quad By \rf{C1}{Lemma 4.51}, we know that the single birth process is recurrent iff the equation
\begin{equation}\label{e3.6}
x_i=\sum_{k\ne 0}\Pi_{ik}x_k,\qquad 0\le x_i\le 1,\qquad i\ge 0
\end{equation}
has only zero solution, where $\Pi_{ik}=(1-\delta_{ik})q_{ik}/q_i$. It is easily seen that equation (\ref{e3.6}) has a non-trivial solution iff the equation
$$x_i=\sum_{k\ne 0}\Pi_{ik}x_k,\qquad i\ge 0; \qquad x_0=1$$
has a nonnegative bounded solution. The following fact will be used several times below:
\begin{equation}\label{e3.6-1}
x_i=\sum_{k\ne i, i_0}\frac{q_{ik}}{q_i-\lambda}x_k +\frac{\gamma_i}{q_i-\lambda}
\Longleftrightarrow (Qx)_i+\lambda x_i=q_{ii_0}(1-\delta_{ii_0})x_{i_0}-\gamma_i,\end{equation}
where $\lambda\in {\mathbb R}$ satisfying some suitable condition. Certainly, here we preassume that
$x_i\in {\mathbb R}$ for every $i\in E$.
By using this fact with $\lambda=0$ and $i_0=0$, we can rewrite the previous equation as
$$(Qx)_0=0,\qquad (Qx)_i=q_{i0},\qquad i\ge 1;\qquad x_0=1.$$
Applying Theorem \ref{t2.1} to $c_i\equiv 0$ and $f_i=q_{i0}(1-\delta_{i0})$, we obtain
the unique solution as follows
$$x_0=1,\quad x_n=1+\sum_{k=1}^{n-1}\sum_{j=1}^k\frac{F_k^{(j)}q_{j0}}{q_{j,j+1}}
=1+\sum_{k=1}^{n-1}\sum_{j=1}^k\frac{F_k^{(j)}q_j^{(0)}}{q_{j,j+1}},\qquad n\ge 1.$$
By (\ref{e0.4-0}), it follows that
$$x_n=1+\sum_{k=1}^{n-1}F_k^{(0)}=\sum_{k=0}^{n-1}F_k^{(0)},\qquad n\ge 1.$$
Clearly, $(x_n)$ is bounded iff $\sum_{k=0}^\infty F_k^{(0)}<\infty$. In other words,
equation (\ref{e3.6}) has only a trivial solution iff $\sum_{k=0}^\infty F_k^{(0)}=\infty$.
The assertion is now proven.\hfill$\square$
\\[0.1cm]

\noindent{\bbb Extinction/return probability}\\[0.1cm]

\noindent
For the remainder of this section, we study the extinction probability.
Here the extinction time $\tau_0$ is the first hitting time of the state $0$. Thus, this
topic is actually a refinement of what studied in the last proposition, in which we pay attention
only on the result either ${\mathbb P}_n[\tau_0<\infty]=1$ or $<1$ rather than its distribution. We will come back this point
after the proof of the next proposition.
For the extinction problem, the rates $q_{0j}\,(j\ne 0)$ play no rule, so one may assume
the state $0$ to be an absorbing state. In other words, we may reduce the state space from $E$ to
$E_1:=\{1, 2, \ldots\}$, and regard the rate $q_{i0}\,(i\ne 0)$ as a killing from $i$. Then we need to
redefine the sequences $\big(\tilde q_{n}^{(k)}\big)$ and $\big(\widetilde F_{n}^{(k)}\big)$ starting
from $1$ but not $0$. However, for our convenience, we prefer to keep the notation $E$,
$\big(\tilde q_{n}^{(k)}\big)$, $\big(\widetilde F_{n}^{(k)}\big)$ and so on. For this, it is better to use the return
time $\sigma_0$ instead of the hitting time $\tau_0$. In the case that the state $0$ is really an absorbing one, we can add a positive rate $q_{01}$ and assume
 that the enlarged process becomes irreducible. Then, the solution of
${\mathbb P}_n[\sigma_0<\infty]$ restricted on $E_1$ gives us the answer of
${\mathbb P}_n[\tau_0<\infty]$ on $E_1$ (as a trivial application of the localization theorem \rf{HG}{Theorem 3.4.1} or \rf{C1}{Theorem 2.13}), so we can return to our original problem.

We remark that in the context of denumerable Markov processes, the topic of this section and much more problems were well
studied in \rf{HG}{Chapter IX}. In the present special case, for the single birth processes,
the problem was studied in \rf{A91}{Chapter 9} or \cite{B86}, using a different technique.

\begin{prop}{  Let the single birth $Q$-matrix $Q=(q_{ij})$ be non-explosive and irreducible. Then the return/extinction probability is as follows:
$$\mathbb P_0(\sigma_0<\infty)=\frac{\sum_{k=1}^{\infty}F_k^{(0)}}{\sum_{k=0}^\infty F_k^{(0)}},\qquad \mathbb P_n(\sigma_0<\infty)=\frac{\sum_{k=n}^{\infty}F_k^{(0)}}{\sum_{k=0}^\infty F_k^{(0)}},\qquad n\ge 1.$$
Furthermore, $\mathbb P_n(\sigma_0<\infty)=1$ for all $n\ge 0$ iff $\mathbb P_0(\sigma_0<\infty)=1$, equivalently iff $\sum_{n=0}^\infty F_n^{(0)}=\infty$.
}\end{prop}

\noindent{\it Proof}\quad By \rf{C1}{Lemma 4.46} with $H=\{0\}$,
 $(\mathbb P_i(\sigma_0<\infty): i\in E)$ is the minimal nonnegative solution to the equation
$$x_i=\sum_{k\ne 0, i}\frac{q_{ik}}{q_i}x_k + \frac{q_{i0}}{q_i}(1-\delta_{i0}),\qquad i\in E.$$
The study on recurrence usually starts from here, the lemma \rf{C1}{Lemma 4.51} used in the last proof simplifies
our study on the recurrence problem, as we have just seen above. By (\ref{e3.6-1}), the last equation is equivalent to
$$(Qx)_i=q_{i0}(1-\delta_{i0})(x_0-1),\qquad i\ge 0.$$
Applying Theorem \ref{t2.1} to $c_i\equiv 0$ and $f_i= q_{i0}(1-\delta_{i0})(x_0-1)$,
we obtain the solution to the last equation:
$$ \aligned
x_n&=x_0+\sum_{0\le k\le n-1}\sum_{0\le j\le k}\frac{F_k^{(j)}}{q_{j,j+1}}q_{j0}(1-\delta_{j0})(x_0-1)\\
&=x_0\bigg\{1+\sum_{1\le k\le n-1}\sum_{1\le j\le k}\frac{F_k^{(j)}}{q_{j,j+1}}q_{j}^{(0)}\bigg\}
   - \sum_{1\le k\le n-1}\sum_{1\le j\le k}\frac{F_k^{(j)}}{q_{j,j+1}}q_{j}^{(0)}\\
&=x_0 \bigg(1+\sum_{1\le k\le n-1} F_k^{(0)}\bigg)-\sum_{1\le k\le n-1}F_k^{(0)},\qquad n\ge 0\quad\text{\rm (by (\ref{e0.4-0}))}.
\endaligned$$
Because $x_n>0$, it follows that
$$x_0\ge \sup_{n\ge 1}\frac{\sum_{k=1}^{n-1}F_k^{(0)}}{\sum_{k=0}^{n-1}F_k^{(0)}}=\sup_{n\ge 1}\frac{\sum_{k=0}^{n-1}F_k^{(0)}-1}{\sum_{k=0}^{n-1}F_k^{(0)}}
=1-\frac 1{\sum_{k=0}^\infty F_k^{(0)}}.$$
From here, we obtain the minimal nonnegative solution:
$$x_0^*=1-\frac 1{\sum_{k=0}^\infty F_k^{(0)}},\qquad x_n^*=1-\frac{\sum_{k=0}^{n-1}F_k^{(0)}}{\sum_{k=0}^\infty F_k^{(0)}},\qquad n\ge 1.$$
We have thus proved the first assertion. The second one is obvious.
\hfill$\square$

\vskip0.2cm
Rewrite the solution just obtained as follows.
$$1-x_0^*=\frac 1{\sum_{k=0}^\infty F_k^{(0)}},\qquad 1-x_n^*=\frac{\sum_{k=0}^{n-1}F_k^{(0)}}{\sum_{k=0}^\infty F_k^{(0)}},\qquad n\ge 1.$$
Renormalize them so that the initial value becomes $1$:
$$x_0=1,\qquad x_n=\sum_{k=0}^{n-1}F_k^{(0)},\qquad n\ge 1$$
which is what we obtained in the last proof. We have thus seen
the relation between the last two propositions.

The study on the Laplace transform of extinction/return time is delayed to
Section 7 (Proposition \ref{pl} which is based on
Lemma \ref{l01}).
\\[4mm]

\setcounter{equation}{0}
\setcounter{section}{5}
\setcounter{thm}{0}
\noindent{\bbb 5\quad Ergodicity, strong ergodicity, and the first moment of return time}\\[0.1cm]

\noindent
Let $E={\mathbb Z}_+$ and $H\subset E$, $H\ne\emptyset, E$. Define $\sigma_H=\inf\{t\ge \eta_1: X(t)\in H\}$,
where $\eta_1$ is the first jump of the process.
When $H$ is a singleton, $H=\{0\}$, for instance, denote $\sigma_{\{0\}}$ by $\sigma_0$ for simplicity.
We now consider the first moment of  the return time $\sigma_0$. To do so, we introduce the following lemma
(cf. \rf{HG}{Lemma 9.4.1}).

\begin{lem}\label{l001}{  Let $(q_{ij})$ be irreducible and assume that its $Q$-process is recurrent.
Then $(x_i^*:= {\mathbb E}_i \sigma_H: i\in E)$ is the minimal nonnegative solution
(may be infinite) to the equation
$$x_i=\frac{1}{q_i}\sum_{k\notin H\cup \{i\}}
  q_{ik}x_k+\frac{1}{q_i}, \qquad i\in E, $$
where $1\cdot\infty=\infty$ and $0\cdot\infty=0$ by convention.
}\end{lem}

\noindent{\it Proof}\quad Let $(y_i^*: i\in E)$ be the minimal nonnegative solution to the equation
$$y_i=\frac{1}{q_i}\sum_{k\notin H\cup \{i\}}
  q_{ik}y_k+\frac{1}{q_i}, \qquad i\in E. $$
By assumption and \rf{C1}{Lemma 4.46}, the quantity $f_{iH}$ defined there is equal
to $1$ for every $i\in E$. Then, $(y_i^*: i\in E)$ coincides with $(e_{iH}(0): i\in E)$
used in \rf{C1}{Lemma 4.48}. Note that $e_{iH}(0)=\int_0^\infty\mathbb P_i(\sigma_H>t)\,\text{d} t
=\mathbb E_i\sigma_H$. The assertion now follows immediately.\hfill$\square$

\vskip0.2cm
In what follows, we use often another sequence $({\tilde d}_n)$ similar to $\big({\widetilde m}_n\big)$ having different initial value:
\begin{equation} \label{e0.6-0}\tilde d_0=0,\qquad \tilde d_n=\frac 1{q_{n,n+1}}\Bigg(1+\sum_{k=0}^{n-1}\tilde q_n^{(k)}\tilde d_k\Bigg),
\qquad n\ge 1,\end{equation}
where $\tilde q_n^{(k)}$ is defined in (\ref{ee0.02}). By Corollary \ref{c0-1}, we have
\begin{equation} \label{e0.6-1}
\tilde d_n=\sum_{1\le j\le n}\frac{\widetilde F_n^{(j)}}{q_{j,j+1}},\qquad n\ge 0
\end{equation}
which is very much the same as (\ref{e3.4-0}).
Again, we omit the superscript $\,\tilde{}\,$ everywhere in $(\tilde d_n)$ once $c_i\equiv 0$. Note that if we rewrite
$$\aligned
{\tilde d}_n&=\frac 1{q_{n,n+1}}\Bigg(1+\sum_{1\le k\le n-1}{\tilde q}_n^{(k)}\tilde d_k\Bigg), \qquad n\ge 1,\\
{\widetilde F}_n^{(0)}&=\frac 1{q_{n,n+1}}\Bigg({\tilde q}_n^{(0)}+\sum_{1\le k\le n-1}{\tilde q}_n^{(k)}{\widetilde F}_k^{(0)}\Bigg),\qquad n\ge 1,
\endaligned
$$
then it is clear that the sequences $\big({\tilde d}_n\big)_{n\ge 1}$
and $\big({\widetilde F}_n^{(0)}\big)_{n\ge 1}$
are also quite close each other.

The main result in this section is as follows. Refer to \rf{C1}{Theorem 4.52\,(2)},
\rf{A91}{Proposition 2.4}, and \cite{YC, Z1, Z2}.

\begin{prop}\label{p5.2}{  Assume that the single birth $Q$-matrix $Q=(q_{ij})$ is irreducible and corresponding process is recurrent. Then
$$\mathbb E_0\sigma_0=\frac 1{q_{01}}+d,\qquad \mathbb E_n\sigma_0=\sum_{k=0}^{n-1}\big(F_k^{(0)}d-d_k\big),\qquad n\ge 1,$$
where
$$d=\varlimsup_{k\to\infty} \frac{\sum_{n=0}^k d_n}{\sum_{n=0}^k F_n^{(0)}}
=\lim_{n\to\infty} \frac{d_n}{F_n^{(0)}}\;\;\text{\rm if the limit exists}.$$
Furthermore, the process is ergodic $($i.e. positive recurrent$)$ iff $d<\infty$;
and it is strongly ergodic iff $\sup_{k\in E} \sum_{n=0}^k\big(F^{(0)}_nd-d_n\big)<\infty$.
Actually, for the last conclusion, the recurrence assumption can be replaced by the
uniqueness one.
}\end{prop}

\noindent{\it Proof}\quad Let $H=\{0\}$. By Lemma \ref{l001}, $({\mathbb E}_i \sigma_0: i\in E)$ is
the minimal nonnegative solution $(x_i^*)$ to the equation
\begin{equation}\label{e6.0}
x_i=\frac{1}{q_i}\sum_{k\notin  \{0, i\}}
  q_{ik}x_k+\frac{1}{q_i}, \qquad i\in E. \end{equation}
Suppose for a moment that $x_i^*<\infty$ first for some
$i\in E$ and then for all $i$ by irreducibility. Next, let
$(x_i)$ be a (finite) solution to (\ref{e6.0}). Then, by (\ref{e3.6-1}), we have
$$(Qx)_i=q_{i0}x_0-1,\qquad i\ge 1;\qquad (Qx)_0=-1.$$
Applying Theorem \ref{t2.1} to $c= 0$ and $f_i=q_{i0}(1-\delta_{i0})x_0-1$ ($i\ge 0$),
we obtain the solution to the last equation:
$$x_n=x_0+\sum_{k=0}^{n-1}\sum_{j=0}^k\frac{F_k^{(j)}f_j}{q_{j,j+1}}=x_0\Bigg(1+\sum_{k=1}^{n-1}\sum_{j=1}^k\frac{F_k^{(j)}q_{j0}}{q_{j,j+1}}\Bigg)-\sum_{k=0}^{n-1}\sum_{j=0}^k\frac{F_k^{(j)}}{q_{j,j+1}},\qquad n\ge 1.$$
By (\ref{e0.4-0}) and (\ref{e0.6-1}), we obtain
$$x_n=x_0\sum_{k=0}^{n-1}F_k^{(0)}-\sum_{k=0}^{n-1}\bigg(\frac {F_k^{(0)}}{q_{01}}+d_k\bigg)
=\sum_{k=0}^{n-1}\bigg[F_k^{(0)}\bigg(x_0-\frac{1}{q_{01}}\bigg)-d_k\bigg],\qquad n\ge 1.$$
Since $x_n>0$, it follows that
$$x_0\sum_{k=0}^{n-1}F_k^{(0)}> \sum_{k=0}^{n-1}\bigg(\frac {F_k^{(0)}}{q_{01}}+d_k\bigg),\qquad n\ge 1.$$
This gives us
$$x_0\ge \sup_{n\ge 1}\frac{\sum_{k=0}^{n-1}({F_k^{(0)}}/{q_{01}}+d_k)}{\sum_{k=0}^{n-1}F_k^{(0)}}
=\frac 1{q_{01}}+\sup_{n\ge 1}\frac{\sum_{k=0}^{n-1} d_k}{\sum_{k=0}^{n-1}F_k^{(0)}}.$$
Now, the minimal property implies that
$$x_0^*
=\frac 1{q_{01}}+\sup_{n\ge 1}\frac{\sum_{k=0}^{n-1} d_k}{\sum_{k=0}^{n-1}F_k^{(0)}}$$
and then
$$x_n^*=\sum_{k=0}^{n-1}\left(F_k^{(0)}\sup_{n\ge 1}\frac{\sum_{j=0}^{n-1} d_j}{\sum_{j=0}^{n-1}F_j^{(0)}}
-d_k\right),\qquad n\ge 1$$
gives us the solution $({\mathbb E}_i \sigma_0: i\in E)$. We claim that the supremum in the last
line has to achieved at infinity. Otherwise, if it is achieved at some finite $n_0$:
$$\frac{\sum_{j=0}^{n_0-1} d_j}{\sum_{j=0}^{n_0-1}F_j^{(0)}}
=\sup_{n\ge 1}\frac{\sum_{j=0}^{n-1} d_j}{\sum_{j=0}^{n-1}F_j^{(0)}}.$$
Then
$$x_0^*
=\frac 1{q_{01}}+ \frac{\sum_{j=0}^{n_0-1} d_j}{\sum_{j=0}^{n_0-1}F_j^{(0)}}$$
and furthermore, $x_{n_0}^*=0$ which is a contradiction with $x_i^*={\mathbb E}_i \sigma_0>0$. Therefore,
$$\sup_{n\ge 1}\frac{\sum_{j=0}^{n-1} d_j}{\sum_{j=0}^{n-1}F_j^{(0)}}
=\varlimsup_{n\to\infty} \frac{\sum_{j=0}^{n} d_j}{\sum_{j=0}^{n}F_j^{(0)}}=: d$$
as required. The next limit in the expression of $d$ is an application of Stolz's Theorem. Now $d<\infty$ since
$x_0^*<\infty$ by assumption. To remove the finiteness assumption of $(x_i^*)$, we claim that the expressions in the first assertion for ${\mathbb E}_n\sigma_0(=x_n^*)$ still hold even $x_i^*=\infty$, since then we must have $d=\infty$. If otherwise, $d<\infty$, then by the last assertion of Theorem \ref{t2.1} and (\ref{e3.6-1}), we would obtain a finite solution to (\ref{e6.0}),
which deduces a contradiction to the assumption $x_i^*=\infty$
by the comparison theorem for the nonnegative solutions
(cf. \rf{C1}{Theorem 2.6}). We have thus proved the first assertion.

Let us remark that the trick used above replacing $\sup_{n\ge 1}$ by $\varlimsup_{n\to\infty}$ was missed in the previous publications. This trick and the one assuming the finiteness of $(x_i^*)$, will be used several times
below but we may not mention it time by time.

Finally, by \rf{C1}{Theorem 4.44}, the single process is ergodic iff $\mathbb E_0\sigma_0<\infty$
which is now equivalent to $d<\infty$. By the same cited theorem, the process is strongly ergodic iff
$\sup_{i\in E}\mathbb E_i\sigma_0<\infty$, equivalently, $\sup_{n\in E} \sum_{k=0}^n\big(F^{(0)}_kd-d_k\big)<\infty$
which follows from the first assertion. As mentioned in the proof of the cited book,
for ergodicity, the uniqueness assumption is enough instead of the recurrence one.
The proof is now finished. \hfill$\square$
\\[4mm]

\setcounter{equation}{0}
\setcounter{section}{6}
\setcounter{thm}{0}
\noindent{\bbb 6\quad Polynomial moments of hitting time and life time}\\[0.1cm]

\noindent{\bbb Polynomial moments of hitting time}\\[0.1cm]

\noindent
We have just studied the first moment of the time of first hitting/return $0$ in the last section.
Now we study the higher-order moments of the first hitting time.

Fix $i_0\ge 0$. Recall that $\sigma_{i_0}$ is the time of first return to $i_0$ after the first jump.
For its higher-moments, we have the following result (cf. \cite{Z3, Z5}).

\begin{prop}{  Assume that the single birth $Q$-matrix $Q=(q_{ij})$ is irreducible and the
corresponding process is $(\ell-1)$-ergodic ($\ell\ge 1$), i.e. ${\mathbb E}_i\sigma_{i_0}^{\ell-1}<\infty$ for every $i\ge 0$. When $\ell=1$, assume additionally that
the process is unique. Then we have
$$\mathbb E_n\sigma_{i_0}^\ell=\begin{cases}\ell \sum_{n\le k\le i_0-1} v_k^{(\ell)}+\big[1-\sum_{n\le k\le i_0-1} u_k\big]{\mathbb E}_{i_0}\sigma_{i_0}^\ell,\  & 0\le n\le i_0;\\
-\ell\sum_{i_0\le k\le n-1} v_k^{(\ell)}+\big[1+\sum_{i_0\le k\le n-1}u_k\big]{\mathbb E}_{i_0}\sigma_{i_0}^\ell,\  & n> i_0;\end{cases}$$
where
\begin{gather}
u_k={\begin{cases}
\sum_{j=i_0-1}^k {q_{j, j+1}}^{-1} {F_k^{(j)}} q_{ji_0}(1-\delta_{ji_0}),\qquad &k\ge i_0,\nonumber\\
1, &k=i_0-1,\\
0,  &0\le k\le i_0-2\end{cases}}\\
v_k^{(\ell)}=\sum_{j=0}^k \frac{F_k^{(j)}}{q_{j, j+1}} {\mathbb E}_j\sigma_{i_0}^{\ell-1},\qquad k\ge 0,\nonumber
\end{gather}
$$\aligned
{\mathbb E}_{i_0}\sigma_{i_0}^\ell&=\ell\varlimsup_{n\to\infty} \bigg({\sum_{i_0\le k\le n} v_k^{(\ell)}}\bigg)\bigg[{1+\sum_{i_0\le k\le n}u_k}\bigg]^{-1}\\
&=\ell\lim_{n\to\infty} \frac{v_n^{(\ell)}}{u_n}\;\; \text{\rm if the limit exists}.\endaligned$$
}\end{prop}

\noindent{\it Proof}\quad By \rf{HG}{Theorem 9.3.3} (cf. \rf{C1}{Proposition 4.56}, or \rf{M1}{Theorem 3.1}),
$(y_i^*:={\mathbb E}_i\sigma_{i_0}^\ell: i\in E)$ is the the minimal nonnegative solution to the following equation:
$$y_i=\sum_{k\ne i, i_0}\frac{1}{q_i} q_{ik}y_k+\frac{\ell}{q_i} {\mathbb E}_i \sigma_{i_0}^{\ell-1},\qquad i\in E.$$
As remarked in the last section, we may assume that $y_i^*<\infty$ for every $i\in E$.
Then, by (\ref{e3.6-1}), we obtain the Poisson equation:
$$(Qy)_i= q_{ii_0}(1-\delta_{ii_0})y_{i_0}-\ell\mathbb E_i\sigma_{i_0}^{\ell-1},\qquad i\in E. $$
Applying Theorem \ref{t2.1} to $c= 0$ and
$f_i= q_{ii_0}(1-\delta_{ii_0})y_{i_0}-\ell\mathbb E_i\sigma_{i_0}^{\ell-1}$, it follows that
 the solution to the last equation is as follows:
$$y_n=y_0+\sum_{0\le k\le n-1}\sum_{j=0}^k\frac{F_k^{(j)}f_j}{q_{j,j+1}}
=y_0+ y_{i_0}\!\!\sum_{0\le k\le n-1}u_k-\ell\!\! \sum_{0\le k\le n-1} v_k^{(\ell)},\qquad n\ge 0.$$
Here in the summation of $u_k$, we have used the character of single birth:
$q_{ji_0}(1-\delta_{ji_0})>0$ only if either $j=i_0-1$ or $j\ge i_0+1$.
In particular, by setting $n=i_0$, it follows that
$$y_0=\ell \sum_{0\le k\le i_0-1} v_k^{(\ell)}+ y_{i_0} \bigg(1-\sum_{0\le k\le i_0-1} u_k\bigg).$$
Return to the original $y_n$, we get
\begin{align}\label{e0.3-1}
y_n&=\ell \bigg[\sum_{0\le k\le i_0-1} v_k^{(\ell)}- \sum_{0\le k\le n-1} v_k^{(\ell)}\bigg]
+y_{i_0} \bigg[1-\sum_{0\le k\le i_0-1} u_k+\sum_{0\le k\le n-1}u_k\bigg]\nonumber\\
&={\begin{cases}
-\ell\sum_{i_0\le k\le n-1} v_k^{(\ell)}+y_{i_0} \big[1+\sum_{i_0\le k\le n-1}u_k\big],\quad & n\ge i_0+1\\
\ell \sum_{n\le k\le i_0-1} v_k^{(\ell)}+y_{i_0} \big[1-\sum_{n\le k\le i_0-1} u_k\big], \quad & n\le i_0.
\end{cases}} \end{align}
When $n\le i_0$, since $\sum_{k\le i_0-1} u_k\le 1$ by definition of $(u_k)$, it is clear that $y_n>0$. When $n\ge i_0+1$, for $y_n>0$, one requires the condition
$$y_{i_0}> \frac{\ell\sum_{i_0\le k\le n-1} v_k^{(\ell)}}{1+\sum_{i_0\le k\le n-1}u_k}$$
and then
$$y_{i_0}\ge \sup_{n\ge i_0+1} \frac{\ell\sum_{i_0\le k\le n-1} v_k^{(\ell)}}{1+\sum_{i_0\le k\le n-1}u_k}.$$
By a reason explained in the last section, this leads to
$$y_{i_0}^*= \ell\varlimsup_{n\to\infty} \frac{\sum_{i_0\le k\le n} v_k^{(\ell)}}{1+\sum_{i_0\le k\le n}u_k}$$
which gives us ${\mathbb E}_{i_0}\sigma_{i_0}^\ell$. Combining it with (\ref{e0.3-1}), we obtain the required assertion.
The limit in ${\mathbb E}_{i_0}\sigma_{i_0}^\ell$ is again an application of Stolz's Theorem since
$\sum_k u_k=\infty$ by the recurrence of the process. To see the last assertion, define a single birth
process on $\{i_0, i_0+1, \ldots\}$ (regarding the set $\{0, 1, \ldots, i_0\}$ as a single state) with rates
$${\bar q}_{ij}=
\begin{cases}
q_{ij} \quad &\text{if } j\ge i_0+1\\
\sum_{k\le i_0}q_{ik} &\text{if } j= i_0,\;\;\;\;\;\; i\ge i_0.
\end{cases}$$
Then $(\bar q_{ij})$ is irreducible and recurrent because so is $(q_{ij})$. Next,
as in (\ref{ee0.01}), we can define a sequence $\big({\widebar F}_k^{(j)}\big)$ on $\{i_0, i_0+1, \ldots\}$.
By induction, it is easy to check that ${\widebar F}_k^{(j)}= {\widetilde F}_k^{(j)}$ for every $k\ge j\ge i_0$.
Hence we have
$$\sum_k{\widebar F}_k^{(i_0)}=\sum_k {\widetilde F}_k^{(i_0)}=\infty$$
by Proposition \ref{pr}. It should be now easy to see that $\sum_k u_k=\infty$ as claimed.
\hfill$\square$
\\[0.1cm]

\noindent{\bbb Polynomial moments of life time}\\[0.1cm]

\noindent
Recall that $\tau_n$ is the time of first hitting the state $n$. If we start from $i\le n-1$, then $\tau_n$ coincides
with the time of fist hitting the set $\{n, n+1, \ldots\}$. For the remainder of this section,
we are going to study the time $\tau_\infty:=\lim_{n\to\infty}\tau_n$.
Next, because $\tau_\infty$ is actually equal to the life time $\eta:=\lim_{n\to\infty}\eta_n$ almost everywhere, where $\{\eta_n\}$ are the successive jumping times:
$$\eta_0\equiv 0,\qquad \eta_n=\inf\{t\ge \eta_{n-1}: X(t)\ne X(\eta_{n-1})\},\qquad n\ge 1,$$
therefore, $\tau_\infty=\infty$ a.e. if
the single birth $Q$-matrix is non-explosive. Thus, the study on the moments of $\tau_\infty$ is
meaningful only for explosive single birth $Q$-matrix. The next result is taken from \cite{Z5}.

\begin{prop}{  Let the single birth $Q$-matrix $Q=(q_{ij})$ be irreducible and explosive
$($i.e. $\sum_n m_n<\infty$ by Proposition $\ref{pr3-1})$. Assume that
the minimal process has finite $(\ell-1)$-th moments of $\tau_\infty$ for some integer $\ell\ge 1$
$($i.e. $E_i\tau_\infty^{\ell-1}<\infty$ for all $i\ge 0)$. Then
$$\mathbb E_n\tau_\infty^\ell=\ell\sum_{k\ge n}\widebar m_k^{(\ell)},\qquad n\ge 0,$$
where
$$\widebar m_n^{(\ell)}=\frac{1}{q_{n,n+1}}\bigg[{\mathbb E}_n\tau_{\infty}^{\ell-1}+\sum_{0\le k\le n-1}q_n^{(k)}{\widebar m}_k^{(\ell)}\bigg]
=\sum_{j=0}^n\frac{F_n^{(j)}\mathbb E_j\tau_\infty^{\ell-1}}{q_{j,j+1}},\qquad n\ge 0.$$
}\end{prop}

\noindent{\it Proof}\quad The last equality of $\widebar m_n^{(\ell)}$ comes from Corollary \ref{c0-1}.
By \rf{C1}{Proposition 4.56} or \cite{M2}, we know that $(\mathbb E_i\tau_\infty^\ell: i\in E)$
is the the minimal nonnegative solution $(y_i^*: i\in E)$ to the following equation:
$$y_i=\sum_{k\ne i}\frac{1}{q_i} q_{ik}y_k+\frac{\ell}{q_i} {\mathbb E}_i \tau_\infty^{\ell-1},\qquad i\in E.$$
That is,
$$(Qy)_i= -\ell\mathbb E_i\tau_\infty^{\ell-1},\qquad i\in E. $$
Applying Theorem \ref{t2.1} to $c= 0$ and $f_i=-\ell\mathbb E_i\tau_\infty^{\ell-1}$ ($i\ge 0$), it follows that
 the solution to the last equation can be expressed as
$$y_n=y_0-\ell\sum_{k=0}^{n-1}\sum_{j=0}^k\frac{F_k^{(j)}\mathbb E_j\tau_\infty^{\ell-1}}{q_{j,j+1}},\qquad n\ge 1.$$
Hence
$$y_n=y_0-\ell\sum_{k=0}^{n-1}\widebar m_k^{(\ell)},\qquad n\ge 1.$$
By the nonnegative and minimal properties, it follows that
$$y_0^*=\sup_{n\ge 1}\Bigg(\ell\sum_{k=0}^{n-1}\widebar m_k^{(\ell)}\Bigg)=\ell\sum_{k=0}^\infty\widebar m_k^{(\ell)}, \qquad y_n^*=\ell\sum_{k=n}^\infty\widebar m_k^{(\ell)},\qquad n\ge 1.$$
Hence, we obtain
$$\mathbb E_n\tau_\infty^\ell=\ell\sum_{k\ge n}\widebar m_k^{(\ell)},\qquad n\ge 0$$
which is the required assertion.\hfill$\square$
\\[4mm]

\setcounter{equation}{0}
\setcounter{section}{7}
\setcounter{thm}{0}
\noindent{\bbb 7\quad Exponential ergodicity and Laplace transform of return time}\\[0.1cm]

\noindent{\bbb Exponential moments of return time and exponential ergodicity}\\[0.1cm]

\noindent
In this section, we consider the exponential moments of return time. At first, we introduce the following lemma for general $Q$-matrices.

\begin{lem}\label{l01}{  Let $(q_{ij})$ be irreducible and assume that its $Q$-process is recurrent.
Next, let $\lambda\in {\mathbb R}$, $\lambda< q_i$ for every $i\in E$. Then for fixed $H\subset E$, $H\ne\emptyset, E$,
$({\mathbb E}_i \exp(\lambda\sigma_H): i\in E)$ is the minimal
solution to the equation
\begin{equation}\label{el01}
x_i=\frac{1}{q_i-\lambda}\sum_{k\notin H\cup \{i\}}
  q_{ik}x_k+\frac{1}{q_i-\lambda}\sum_{k\in H\setminus \{i\}}
  q_{ik}, \qquad i\in E.
\end{equation}
}\end{lem}

\noindent{\it Proof}\quad Let $(y_i^*: i\in E)$ be the minimal nonnegative solution to the equation
$$y_i=\frac{1}{q_i-\lambda}\sum_{k\notin H\cup \{i\}}
  q_{ik}y_k+\frac{1}{q_i-\lambda}, \qquad i\in E. $$
By the recurrent assumption and \rf{C1}{Lemma 4.46}, the quantity $f_{iH}$ defined there is equal
to $1$ for every $i\in E$. Then, $(y_i^*: i\in E)$ coincides with $(e_{iH}(\lambda): i\in E)$
used in \rf{C1}{Lemma 4.48}. Moreover, by the proof given on \rf{C1}{page 148}, we have
${\mathbb E}_i \exp(\lambda\sigma_H)=1+ \lambda y_i^*$ for every $i\in E$. Besides, it can be checked that $(1+ \lambda y_i^*:i\in E)$ is a nonnegative solution to equation (\ref{el01}).
Hence ${\mathbb E}_i \exp(\lambda\sigma_H)=1+ \lambda y_i^*\ge x_i^*$ for every $i\in E$, where  $(x_i^*: i\in E)$ is the minimal nonnegative solution to equation (\ref{el01}). We are now going to prove that ${\mathbb E}_i \exp(\lambda\sigma_H)=x_i^*$ for all $i\in E$.
The proof is split into two parts: either $\lambda\ge 0$ or $\lambda< 0$.

First, let $\lambda\ge 0$. It is easily seen that
$(x_i^*-1:i\in E)$ is a nonnegative solution to the equation
$$y_i=\frac{1}{q_i-\lambda}\sum_{k\notin H\cup \{i\}}
  q_{ik}y_k+\frac{\lambda}{q_i-\lambda}, \qquad i\in E. $$
Hence, $x_i^*-1\ge \lambda y_i^*$ since $(\lambda y_i^*)$ is
the minimal nonnegative solution to the equation above, by the linear combination theorem \rf{C1}{Theorem 2.12\,(1)}. That is, $x_i^*\ge 1+\lambda y_i^*$. Combining what we have proved in the last paragraph, it follows that $x_i^*={\mathbb E}_i \exp(\lambda\sigma_H)$ for all $i\in E$.

Next, let $\lambda< 0$. Denote by $(\bar y_i: i\in E)$ the minimal nonnegative solution to the equation
\begin{equation}\label{el01-1}
y_i=\frac{1}{q_i-\lambda}\sum_{k\notin H\cup \{i\}}
  q_{ik}y_k+\bigg[1- \frac{1}{q_i-\lambda}\sum_{k\notin H\cup \{i\}}q_{ik}\bigg], \qquad i\in E.
\end{equation}
Clearly, we have $\bar y_i\le 1$ since $y_i\equiv 1$ is a solution to the equation.
We claim that $\bar y_i\equiv 1$. To see this, note that
$(1-\bar y_i: i\in E)$ is the maximal solution to the equation
\begin{equation}\label{el01-2}
y_i=\frac{1}{q_i-\lambda}\sum_{k\notin H\cup \{i\}}q_{ik}y_k,\qquad 0\le y_i\le 1,\qquad i\in E.\end{equation}
By a comparison lemma \rf{C1}{Lemma 3.14}, it suffices to show that
the equation
$$y_i=\frac{1}{q_i}\sum_{k\notin H\cup \{i\}}q_{ik}y_k,\qquad 0\le y_i\le 1,\qquad i\in E$$
has only trivial (i.e. zero-) solution. Then this follows by
the recurrence assumption and \rf{C1}{Lemma 4.46}. We remark that there is an alternative way
to prove that $\bar y_i\equiv 1$, using the uniqueness rather than the recurrence assumption.
Actually, equation (\ref{el01-2}) is an exit equation for a modified $Q$-matrix (any local modification of a $Q$-matrix does not interfere the uniqueness). The exit solution to
(\ref{el01-2}) should be zero by uniqueness assumption.

We now return to our main proof.  By the linear combination theorem \rf{C1}{Theorem 2.12\,(1)},
$(x_i^*-\lambda y_i^*:i\in E)$ is the minimal nonnegative solution to equation (\ref{el01-1}).
Hence $x_i^*-\lambda y_i^*=\bar y_i\equiv 1$ as we have just proved in the last paragraph.
Therefore we conclude that $x_i^*=1+\lambda y_i^*={\mathbb E}_i \exp(\lambda\sigma_H)$
for all $i\in E$. We have thus completed the proof of the lemma.
\hfill$\square$

\vskip0.2cm
Now we present our results about the exponential moments of the return time $\sigma_0$, which can be referred in \cite{Z2}.

\begin{prop}\label{pe}{  Let the single birth $Q$-matrix $(q_{ij})$ be irreducible. Assume that its process is ergodic. Define $\big(\widetilde F_k^{(i)}\big)$ and $\big(\tilde d_k\big)$ by setting $c_i\equiv\lambda>0$. Then for small $\lambda$,
$$\mathbb E_0\text{e}^{\lambda\sigma_0}\!=\!\frac {q_{01} (1+\lambda\tilde d)}{q_{01}-\lambda}\!<\!\infty\;\;\text{  and}\;\;
\mathbb E_n\text{e}^{\lambda\sigma_0}\!=\!1+ \lambda \sum_{k=0}^{n-1}\Big(\widetilde F_k^{(0)}\tilde d-\tilde d_k\Big)\!<\!\infty,\quad n\ge 1$$
iff
$${\tilde d}:=\varlimsup_{n\to \infty}
\mathbbm{1}_{\big\{{\sum_{k=0}^{n}\widetilde F_k^{(0)}}>0\big\}}
\frac{\sum_{k=0}^{n}\tilde d_k}{\sum_{k=0}^{n}\widetilde F_k^{(0)}}<\infty$$
and
\begin{equation}\label{e4.4-0}
{\tilde d}\, \sum_{k=0}^{n-1}\widetilde F_k^{(0)}>\sum_{k=0}^{n-1} \tilde d_k\quad \text{  whenever } \sum_{k=0}^{n-1}\widetilde F_k^{(0)}\le 0 \quad\text{  for } n\ge 2.
\end{equation}
Furthermore, once ${\widetilde F}_n^{(0)}>0$ for large enough $n$ and $\sum_n {\widetilde F}_n^{(0)}=\infty$, we have $${\tilde d}=\lim_{n\to\infty} \frac{{\tilde d}_n}{\widetilde F_n^{(0)}}\;\;\text{\rm if the limit exists}.$$
Finally, the process is exponentially ergodic iff both $\tilde d<\infty$ and $(\ref{e4.4-0})$ holds.
}\end{prop}

\noindent{\it Proof}\quad Let $\lambda\in (0, q_i)$ for every $i\in E$ and set $H=\{0\}$.
Then by Lemma \ref{l01}, $(\mathbb E_i\text{e}^{\lambda\sigma_0}: i\in E)$ is the minimal solution $(x_i^*)$ of the following equation
$$x_i=\frac{1}{q_i-\lambda}\sum_{k\notin \{0,i\}}
  q_{ik}x_k+\frac{q_{i0}(1-\delta_{i0})}{q_i-\lambda}, \qquad x_i\ge1,\qquad i\in E. $$
Assume that $x_i^*<\infty$ for every $i\in E$ for a moment, and let
$(x_i)$ be a finite nonnegative solution to the last equation. Then, by (\ref{e3.6-1}), we have
\begin{equation}\label{e4.4}
(Qx)_i+\lambda x_i=q_{i0}(x_0-1),\qquad i\ge 1;\qquad (Qx)_0+\lambda x_0=0.
\end{equation}
Applying Theorem \ref {t2.1} to $c_i\equiv\lambda$ and $f_i=q_{i0}(1-\delta_{i0})(x_0-1)$ for all $i\ge 0$, we obtain
$$\aligned x_n&=x_0\Bigg(1-\lambda\sum_{k=0}^{n-1}\sum_{j=0}^k\frac{\widetilde F_k^{(j)}}{q_{j,j+1}}\Bigg)+(x_0-1)\sum_{k=1}^{n-1}\sum_{j=1}^{k}\frac{\widetilde F_k^{(j)}q_{j0}}{q_{j,j+1}}\\
&=x_0\Bigg(1-\lambda\sum_{k=0}^{n-1}\sum_{j=0}^k\frac{\widetilde F_k^{(j)}}{q_{j,j+1}}\Bigg)+(x_0-1)\sum_{k=1}^{n-1}\sum_{j=1}^{k}\frac{\widetilde F_k^{(j)}(\tilde q_{j}^{(0)}+\lambda )}{q_{j,j+1}},\qquad n\ge 1.\endaligned$$
Due to the explicit representation of ${\widetilde F}_n^{(k)}$, ${\widetilde m}_n$ and $\tilde d_n$,
given in (\ref{e0.4-0}), (\ref{e3.4-0}) and (\ref{e0.6-1}) respectively, we have not only
\begin{equation}\label{e4.4-1}\widetilde m_n=\sum_{0\le j\le n}\frac{\widetilde F_n^{(j)}}{q_{j,j+1}}=\frac 1{q_{01}}\widetilde F_n^{(0)}+\tilde d_n,\qquad n\ge 0\end{equation}
but also that
\begin{align}\label{e4.4-3} x_n&=x_0\Bigg(1-\lambda\sum_{k=0}^{n-1}\widetilde m_k\Bigg)+(x_0-1)\sum_{k=1}^{n-1}\big(\widetilde F_k^{(0)}+\lambda\tilde d_k\big)\nonumber\\
&=x_0\bigg(1-\frac \lambda{q_{01}}\bigg)\sum_{k=0}^{n-1}\widetilde F_k^{(0)}-\sum_{k=0}^{n-1}\big(\widetilde F_k^{(0)}+\lambda\tilde d_k\big)+1,\qquad n\ge 1.
\end{align}
Since $x_n> 1$, we get
$$x_0\bigg(1-\frac \lambda{q_{01}}\bigg)\sum_{k=0}^{n-1}\widetilde F_k^{(0)}> \sum_{k=0}^{n-1}\big(\widetilde F_k^{(0)}+\lambda\tilde d_k\big),\qquad n\ge 1.$$
That is
\begin{equation}\label{e4.4-2}\bigg[x_0\bigg(\frac{1}{\lambda}-\frac{1}{q_{01}}\bigg)-\frac{1}{\lambda}\bigg]
\sum_{k=0}^{n-1}\widetilde F_k^{(0)}>\sum_{k=0}^{n-1} \tilde d_k,\qquad n\ge 1. \end{equation}
Note that on the one hand, if $x_0^*=x_0^*(\lambda_0)<\infty$, then $x_0^*=x_0^*(\lambda)<\infty$ for every $\lambda\in (0, \lambda_0)$,
by the comparison theorem (cf. \rf{C1}{Theorem 2.6}).
On the other hand, when $\lambda=0$, we have
$$\sum_{k=0}^{n}\widetilde F_k^{(0)}=\sum_{k=0}^{n} F_k^{(0)}>0\quad\text{and}\quad
\sum_{k=0}^{n} \tilde d_k =\sum_{k=0}^{n} d_k>0,\qquad n\ge 1.$$
For each fixed $n$, $\sum_{k=0}^{n}\widetilde F_k^{(0)}$ and $\sum_{k=0}^{n} \tilde d_k$ are analytic in $\lambda$,
and so should be positive for sufficient small $\lambda$, say $\lambda\le \lambda_1$ for some $\lambda_1\le\lambda_0$.
Then by (\ref{e4.4-2}), we should have
$$x_0\bigg(\frac{1}{\lambda}-\frac{1}{q_{01}}\bigg)-\frac{1}{\lambda}
>0,\qquad \lambda\in (0, \lambda_1) $$
independent of $n$. Therefore, by the minimal property,  we have
$$x_0^*\bigg(\frac{1}{\lambda}-\frac{1}{q_{01}}\bigg)-\frac{1}{\lambda}
=\varlimsup_{n\to \infty}\mathbbm{1}_{\big\{{\sum_{k=0}^{n}\widetilde F_k^{(0)}}>0\big\}}
\bigg[{\sum_{k=0}^{n}\tilde d_k}\bigg]\bigg[{\sum_{k=0}^{n}\widetilde F_k^{(0)}}\bigg]^{-1}=\tilde d,$$
i.e.
\begin{equation}\label{e4.4-4}\mathbb E_0\text{e}^{\lambda\sigma_0}=x_0^*=\frac{q_{01}(1+\lambda \tilde d)}{q_{01}-\lambda}.\end{equation}
Since $x_0^*$ satisfies (\ref{e4.4-2}), we obtain condition (\ref{e4.4-0}). Then
$$\mathbb E_n\text{e}^{\lambda\sigma_0}=1+ \lambda \sum_{k=0}^{n-1}\Big(\widetilde F_k^{(0)}\tilde d-\tilde d_k\Big),\qquad n\ge 1.$$

Conversely, if $\tilde d<\infty$ and (\ref{e4.4-0}) holds. Then
starting from $x_0=x_0^*$ given in (\ref{e4.4-4}) and defining $x_n$ by (\ref{e4.4-3}),
we obtain a solution $(x_i>1: i\in E)$ to (\ref{e4.4}). By (\ref{e3.6-1}), we obtain a finite nonnegative solution
to the original equation for $\big({\mathbb E}_i e^{\lambda\sigma_0}: i\in E\big)$, and hence
the minimal solution $\big(x_i^*={\mathbb E}_i e^{\lambda\sigma_0}: i\in E\big)$ should be finite.

Finally, by \rf{C1}{Theorem 4.44}, the process is exponentially ergodic iff
$\mathbb E_0\text{e}^{\lambda\sigma_0}$ $<\infty$, equivalently, $\tilde d<\infty$ and (\ref{e4.4-0}) holds.
The last assertion of the proposition then follows.\hfill$\square$

\vskip0.2cm
In contract to the ergodic case, one may study the exponential
decay (in the transient case) for which the Poisson equation
becomes
$$Qg+ \lambda g=0,\qquad g>0.$$
With $c_i\equiv \lambda$, by Theorem \ref{t2.1}, the solution is
$$g_n=g_0\bigg[1-\lambda\sum_{0\le k\le n-1}\sum_{0\le j\le k}\frac{\widetilde F_k^{(j)}}{q_{j,j+1}}\bigg]
=g_0\bigg[1-\lambda\sum_{0\le k\le n-1}{\widetilde m}_k\bigg],
\qquad n\ge 0.$$
This is somehow simpler than the previous one. However, these
two exponential cases are actually much harder than the others,
for instance we do not know at the moment how to remove condition (\ref{e4.4-0}).
That is showing for some $\lambda\!>\!0$, small enough, $\sum_{k=0}^{n}\widetilde F_k^{(0)}\!>\!0$ for all $n$
\big(or equivalently, $\varliminf_{n\to\infty}\sum_{k=0}^{n}\widetilde F_k^{(0)}\!>\!0$\big).
This seems necessary for the exponential ergodicity since
$\sum_{k=0}^{\infty}\widetilde F_k^{(0)}=\infty$ when $\lambda=0$ by the recurrence
(which is much weaker than exponential ergodicity) and
$\lambda$ is allowed to be very small.
Actually, to figure out a criterion, one needs much
more work using different approaches, refer to \rf{C1}{Chapter 9} and \cite{C6}
for some details.
\\[0.1cm]

\noindent{\bbb Laplace transform of the return/extinction time}\\[0.1cm]

\noindent
Note that for negative $\lambda$, $\mathbb E_i\text{e}^{\lambda\sigma_0}$ is the Laplace transform of $\sigma_0$. The proof of Proposition \ref{pe} is still available. So we get the following result.

\begin{prop}\label{pl}{
Define $\big(\widetilde F_k^{(i)}\big)$ and $\big(\tilde d_k\big)$ by $(\ref{ee0.01})$ and $(\ref{e0.6-0})$, respectively, with $c_i\equiv -\lambda<0$.
Let the single birth process be recurrent. Then
the Laplace transform of $\sigma_0$ is given by
$$\mathbb E_0\text{e}^{-\lambda\sigma_0}=\frac {q_{01} (1-\lambda\tilde d)}{q_{01}+\lambda},\qquad \mathbb E_n\text{e}^{-\lambda\sigma_0}=1- \lambda \sum_{k=0}^{n-1}\Big(\widetilde F_k^{(0)}\tilde d-\tilde d_k\Big),\qquad n\ge 1,$$
where
$$\tilde d=\lim_{n\to\infty}\frac{\sum_{k=0}^{n-1}\tilde d_k}{\sum_{k=0}^{n-1}\widetilde F_k^{(0)}}
=\lim_{n\to\infty}\frac{{\tilde d}_n}{ {\widetilde F}_n^{(0)}} \quad \text{\rm if the limit exists}.$$
}\end{prop}

\noindent{\it Proof}\quad Following the proof of Proposition \ref{pe}, replacing $\lambda$ by $-\lambda$, we arrive at
$$\aligned x_n
&=x_0\bigg(1+\frac \lambda{q_{01}}\bigg)\sum_{k=0}^{n-1}\widetilde F_k^{(0)}-\sum_{k=0}^{n-1}\big(\widetilde F_k^{(0)}-\lambda\tilde d_k\big)+1,\\
&=:x_0 \alpha_{n-1}-\beta_{n-1}, \qquad n\ge 1.\endaligned$$
By the minimal nonnegative property,
$x_0^*=\sup_{n\ge 1}\beta_{n}/\alpha_{n}$, and then we indeed have
$$x_0^*=\varlimsup_{n\to\infty}\frac{\beta_{n}}{\alpha_{n}}.$$
We now show that we can replace $\varlimsup_{n\to\infty}$ by $\lim_{n\to\infty}$.
Noting that on the one hand, since $x_n\in (0, 1]$, we have
$$\frac{\beta_{n}}{\alpha_{n}}< x_0\le \frac{\beta_{n}+1}{\alpha_{n}},\qquad n\ge 1.$$
On the other hand, following the proof for
$$\sum_{k}{\widetilde m}_k=\infty\Longleftrightarrow \sum_{k}{m}_k=\infty$$
given in Section 3, we can prove that $\sum_{k}{\widetilde F}_k^{(0)}=\infty$ since $\sum_{k}{F}_k^{(0)}=\infty$ by the recurrent
assumption (i.e. $\gamma_j\equiv 1$).
Hence we can rewrite $\varlimsup_{n\to\infty}{\beta_{n}}/{\alpha_{n}}$
as $\lim_{n\to\infty}{\beta_{n}}/{\alpha_{n}}$.
Therefore, we have
$$\aligned
x_0^*&=\lim_{n\to\infty}\bigg[\sum_{k=0}^{n-1}\big({\widetilde F}_k^{(0)}-\lambda {\tilde d}_k\big) \bigg]\bigg\{\bigg[1+\frac{\lambda}{q_{01}}\bigg]
\sum_{k=0}^{n-1}{\widetilde F}_k^{(0)}\bigg\}^{-1}\\
&= \frac{q_{01}}{q_{01}+\lambda}\lim_{n\to\infty}\bigg[1-\lambda \frac{\sum_{k=0}^{n-1} {\tilde d}_k}{\sum_{k=0}^{n-1}{\widetilde F}_k^{(0)}}\bigg]\\
&= \frac{q_{01}}{q_{01}+\lambda} \big[1-\lambda{\tilde d}\,\big].
\endaligned.$$
Furthermore,
$$x_n^*=(1-\lambda {\tilde d})\sum_{k=0}^{n-1}\widetilde F_k^{(0)}-\sum_{k=0}^{n-1}\big(\widetilde F_k^{(0)}-\lambda\tilde d_k\big)+1
=1- \lambda \sum_{k=0}^{n-1}\Big(\widetilde F_k^{(0)}\tilde d-\tilde d_k\Big),\quad n\ge 1.$$

The last limit in $\tilde d$ is an application of Stolz's Theorem.\hfill$\square$
\\[0.1cm]

\noindent{\bbb Exponential moments and Laplace transform of the life time}\\[0.1cm]

\noindent
Now we return to $\tau_\infty$.

\begin{prop}\label{p7.4}{  Assume that the single birth $Q$-matrix $Q=(q_{ij})$ is explosive and irreducible. Define $(\widetilde m_k)$ by $(\ref{e3.4-00})$ with $c_i\equiv \lambda$.
For the corresponding minimal process,

$(i)$ if there exists a $\lambda>0$ such that $\lambda\sum_{k=0}^{n-1}\widetilde m_k<1$ for every $n>1$,
then
$$\mathbb E_n\text{e}^{\lambda\tau_\infty}=1+ \lambda \Bigg[\bar c\Bigg(1-\lambda\sum_{k=0}^{n-1}\widetilde m_k\Bigg)-\sum_{k=0}^{n-1}\widetilde  m_k\Bigg],\qquad n\ge 0,$$
where
$$\bar c=\varlimsup_{n\to\infty }\frac{\sum_{k=0}^{n}\widetilde m_k}{1-\lambda\sum_{k=0}^{n}\widetilde m_k}.$$
Furthermore, the process decays exponentially fast provided $\bar c<\infty$.

$(ii)$ For $\lambda>0$, the Laplace transform of $\tau_\infty$ is given by
$$\mathbb E_n\text{e}^{-\lambda\tau_\infty}=\frac{1+\lambda\sum_{0\le k\le n-1}\widetilde m_k}{1+\lambda\sum_{k\ge 0}\widetilde  m_k},\qquad n\ge 0.$$
}\end{prop}

\noindent{\it Proof}\quad Define
$$e_{i\infty}(\lambda)=\int_0^\infty\text{e}^{\lambda t}\mathbb P_i(\tau_\infty>t)\,\text{d} t$$
with $\lambda<q_i$ for all $i\ge 0$. Note that the process is explosive and
$$\mathbb E_i\text{e}^{\lambda\tau_\infty}=1+\lambda e_{i\infty}(\lambda).$$
Because $\mathbb P_m(\tau_n<\eta)=1$ for every pair $m<n$, we have $\mathbb P_m(\tau_n<\infty)=1$ and furthermore $\mathbb P_m(\tau_\infty<\infty)=1$ for every $m$, as $n$ goes to $\infty$. Then by \rf{C1}{Lemma 4.48}, $(e_{i\infty}(\lambda))$ is the minimal solution to the equation
$$x_i=\frac{q_i}{q_i-\lambda}\sum_{k}\Pi_{ik}x_k+\frac 1{q_i-\lambda}, \qquad i\ge 0.$$
By (\ref{e3.6-1}), we can rewrite the equation as
$$(Qx)_i+\lambda x_i=-1,\qquad i\ge 0.$$
Applying Theorem \ref{t2.1} to $c_i\equiv \lambda$ and $f_i\equiv -1$, the solution of the equation has
the form:
$$\aligned x_n&=x_0\Bigg(1-\lambda\sum_{k=0}^{n-1}\sum_{j=0}^k\frac{\widetilde F_k^{(j)}}{q_{j,j+1}}\Bigg)
-\sum_{k=0}^{n-1}\sum_{j=0}^k\frac{\widetilde F_k^{(j)}}{q_{j,j+1}}\\
&=x_0\Bigg(1-\lambda\sum_{k=0}^{n-1}\widetilde m_k\Bigg)-\sum_{k=0}^{n-1}\widetilde  m_k,\qquad n\ge 1.\endaligned$$
Note that $\lambda<q_0=q_{01}$ and $\lambda\widetilde m_0<1$. If there exists a positive
$\lambda$ small enough so that $\lambda\sum_{k=0}^{n-1}\widetilde m_k<1$ for every $n>1$,
then by the argument above and the minimal property of the solution, one gets
$$
e_{0\infty}(\lambda)=\sup_{n\ge 1}\frac{\sum_{k=0}^{n-1}\widetilde m_k}{1-\lambda\sum_{k=0}^{n-1}\widetilde m_k}=
\varlimsup_{n\to\infty}\frac{\sum_{k=0}^{n}\widetilde m_k}{1-\lambda\sum_{k=0}^{n}\widetilde m_k}=:\bar c
$$
and
$$e_{n\infty}(\lambda)=\bar c\Bigg(1-\lambda\sum_{k=0}^{n-1}\widetilde m_k\Bigg)-\sum_{k=0}^{n-1}\widetilde  m_k,\qquad n\ge 1.$$
Then the first assertion follows.

For the Laplace transform of $\tau_\infty$, the argument above still works because now we deal with the case of $-\lambda<0$. By the explosive property, we know that
$\sum_{k=0}^\infty\widetilde m_k<\infty$. Hence we have
$$e_{0\infty}(-\lambda)=\bar c=\frac{\sum_{k=0}^\infty\widetilde m_k}{1+\lambda\sum_{k=0}^\infty\widetilde m_k}$$
and
$$e_{n\infty}(-\lambda)=\bar c\Bigg(1+\lambda\sum_{k=0}^{n-1}\widetilde m_k\Bigg)-\sum_{k=0}^{n-1}\widetilde  m_k
=\frac{\sum_{k=n}^\infty \widetilde m_k}{1+\lambda\sum_{k=0}^\infty\widetilde  m_k},\qquad n\ge 1.$$
Finally, we have
$$\mathbb E_n\text{e}^{-\lambda\tau_\infty}=1-\frac{\lambda\sum_{k=n}^\infty \widetilde m_k}{1+\lambda\sum_{k=0}^\infty\widetilde  m_k}
=\frac{1+\lambda\sum_{0\le k\le n-1}\widetilde m_k}{1+\lambda\sum_{k\ge 0}\widetilde  m_k},\qquad n\ge 0.$$
The proof for the second assertion is now finished. \hfill$\square$

\vskip0.2cm
A more careful study on part $(i)$ of Proposition \ref{p7.4}, refer to Proposition \ref{pe}.\\[4mm]

\setcounter{equation}{0}
\setcounter{section}{8}
\setcounter{thm}{0}
\noindent{\bbb 8\quad Examples}\\[0.1cm]

\noindent
In the special case of birth--death processes, the problems studied here have rather complete solutions, see for instance \rf{C1}{Theorem 4.55}. As mentioned in the introduction of the paper, much more models have been studied in the past years. Here we make a little addition. The next example is taken from \cite{BGR}.

\begin{exm}[uniform catastrophes]\label{ex8.1}{\rm
Let
$$q_{i,i+1}=b\, i,\qquad i\ge 0;\qquad q_{ij}=a,\qquad j=0,1,\dots, i-1;$$
and $q_{ij}=0$ for other $j>i+1$, where $a$ and $b$ are positive constants. Then the extinction of the process
has an exponential distribution
$${\mathbb E}_n\text{e}^{-\lambda\tau_0}=\frac{a}{a+\lambda},
\qquad \lambda>0,\; n\ge 1.$$
It is surprising that the distribution is independent of $b$ and the starting point $n$. Redefine $q_{01}=1$. Then the irreducible process is indeed
strongly ergodic.}
\end{exm}

\noindent{\it Proof}\quad We need to consider the case that
$q_{01}>0$ only. With $c_i\equiv -\lambda\in {\mathbb R}$ and then
${\tilde q}_n^{(k)}=(k+1) a +\lambda$ for $k\le n-1$, by using
(\ref{ee0.01}), (\ref{e0.6-0}), and induction, one may check
that
$$\aligned
{\widetilde F}_n^{(0)}&=\frac{a+\lambda}{nb}\prod_{1\le k\le n-1}\Bigg(1+\frac{(k+1)a+\lambda}{kb}\Bigg),
\qquad \prod_{\emptyset}=:1,\\
{\tilde d}_n&= \frac{1}{nb}\prod_{1\le k\le n-1}\Bigg(1+\frac{(k+1)a+\lambda}{kb}\Bigg),\qquad n\ge 1.
\endaligned $$
Since for each fixed $\lambda\in {\mathbb R}$,
$$\log\Bigg(1+\frac{(n+1)a+\lambda}{nb}\Bigg)\to
\log\Bigg(1+\frac{a}{b}\Bigg)>0\qquad\text{as }n\to\infty,$$
we have $\lim_{n\to\infty}{\widetilde F}_n^{(0)}=\infty$ and so
$\sum_n {\widetilde F}_n^{(0)}=\infty$.
As an application of this fact with $\lambda=0$, it follows that the process is recurrent
(Proposition \ref{pr}) and then should be non-explosive
((\ref{e4.4-1}) and Proposition \ref{pr3-1}).

Next, because
$$\sum_n {\widetilde F}_n^{(0)}=\infty,\qquad {\widetilde F}_n^{(0)}=(a+\lambda){\tilde d}_n,\;\; n\ge 1,$$
it follows that
$${\tilde d}=\lim_{n\to\infty} \frac{{\tilde d}_n}{{\widetilde F}_n^{(0)}}=\frac{1}{a+\lambda}.$$
Hence, we have
$${\widetilde F}_n^{(0)}{\tilde d}={\tilde d}_n,\qquad n\ge 1,$$
From here, when $\lambda=0$ in particular, we obtain
$$\sup_k \sum_{n=0}^k \big(F_n^{(0)} d- d_n\big)=d=a^{-1}<\infty. $$
Hence the process is strongly ergodic by Proposition \ref{p5.2}.

By using Proposition \ref{pl}, we obtain
$$\aligned
\mathbb E_0\text{e}^{-\lambda\sigma_0}&=\frac{aq_{01}}{(a+\lambda)(q_{01}+\lambda)},\\
\mathbb E_n\text{e}^{-\lambda\sigma_0}&=1-\lambda {\tilde d}=\frac{a}{a+\lambda}
  =\mathbb E_n\text{e}^{-\lambda\tau_0},\qquad n\ge 1.
\endaligned$$
Therefore, we have proved the first assertion.

Even though it is now automatic that the process
is exponentially ergodic, implied by the strongly ergodicity,
we would like to check the effectiveness of Proposition \ref{pe} for this model. To do so, reset $c_i\equiv \lambda>0$.
Then
$$\aligned
{\widetilde F}_n^{(0)}&=\frac{a-\lambda}{nb}\prod_{1\le k\le n-1}\Bigg(1+\frac{(k+1)a-\lambda}{kb}\Bigg),\\
{\tilde d}_n&= \frac{1}{nb}\prod_{1\le k\le n-1}\Bigg(1+\frac{(k+1)a-\lambda}{kb}\Bigg),\qquad n\ge 1.
\endaligned $$
Clearly, ${\widetilde F}_n^{(0)}>0$ and so does ${\tilde d}_n$ for every $\lambda\in (0, a)$. As we have proved above
$$\sum_n {\widetilde F}_n^{(0)}=\infty,\qquad
{\tilde d}=\lim_{n\to\infty} \frac{{\tilde d}_n}{{\widetilde F}_n^{(0)}}=\frac{1}{a-\lambda}<\infty,$$
and hence the process is exponentially ergodic by
Proposition \ref{pe}. Actually, we have
$$\aligned
\text{\hskip2.5truecm}\mathbb E_0\text{e}^{\lambda\sigma_0}&=\frac{aq_{01}}{(a-\lambda)(q_{01}-\lambda)},\\
\mathbb E_n\text{e}^{\lambda\sigma_0}&=\frac{a}{a-\lambda},
\qquad n\ge 1,\;\;\;\; \lambda\in (0, a\wedge q_{01}).
\text{\hskip2.5truecm}\square\endaligned$$

\begin{exm}\label{ex8.2}{\rm
Consider the single birth $Q$-matrix $(q_{ij})$
with
$$q_{i0}>0, \;\; q_{i, i+1}>0,\;\; q_{ij}=0\;\;\text{\rm for all other }j\ne i.$$
Let $c_i\in {\mathbb R}$. Then
\begin{itemize}
\item [(1)] we have
\begin{align}\label{e8.0}\text{\hskip-1em}
{\widetilde F}_i^{(i)}=1,\quad
&{\widetilde F}_n^{(i)}=\frac{q_{n0}-c_n}{q_{n, n+1}}
\prod_{i+1\le k\le n-1}\bigg[1+ \frac{q_{k0}-c_k}{q_{k, k+1}}\bigg],\\
&\prod_{\emptyset}=:1,\qquad n>i\ge 0,\nonumber\end{align}
and then $(\tilde m_n)$ and $\big(\tilde d_n\big)$ are given by $(\ref{e3.4-0})$ and $(\ref{e0.6-1})$, respectively.
\item [(2)] In particular, if $q_{n0}-c_n\equiv q_{10}-c_1$ for
every $n\ge 1$, then
\begin{gather}
{\widetilde F}_i^{(i)}=1,\;\;
{\widetilde F}_n^{(i)}=\frac{q_{10}-c_1}{q_{n, n+1}}
\prod_{k=i+1}^{n-1}\bigg[1+ \frac{q_{10}-c_1}{q_{k, k+1}}\bigg],\;\prod_{\emptyset}=:1,\;\; n>i\ge 0,\nonumber\\
{\widetilde m}_0=\frac{1}{q_{01}},\;\;{\widetilde m}_n=\frac{1}{q_{n, n+1}}
\prod_{k=0}^{n-1}\bigg[1+ \frac{q_{10}-c_1}{q_{k, k+1}}\bigg],\qquad n\ge 1,\nonumber\\
{\tilde d}_0=0,\;\;
{\tilde d}_n=\frac{1}{q_{n, n+1}}
\prod_{1\le k\le n-1}\bigg[1+ \frac{q_{10}-c_1}{q_{k, k+1}}\bigg],\qquad n\ge 1.\nonumber
\end{gather}
Furthermore, the process is explosive if
$$\kappa':=\lim_{n\to\infty}
\frac{n(q_{n+1, n+2}-q_{n, n+1}-q_{10})}{q_{n, n+1}+q_{10}}>1$$
($q_{n, n+1}=(n+1)^{\gamma}$ for $\gamma> 1$ for example).
Otherwise, if $\kappa'<1$ ($q_{n, n+1}=(n+1)^{\gamma}$ for some $\gamma\le 1$ for instance), then the process is
unique. If so, the process is indeed strongly ergodic.
\end{itemize}
}\end{exm}

\noindent{\it Proof}\quad (a)
By assumption, we have ${\tilde q}_n^{(k)}=q_{n0}-c_n$ for every $k< n$. Hence, by (\ref{ee0.01}), we obtain
\begin{equation}\label{e8.4}{\widetilde F}_n^{(i)}=\frac{{\tilde q}_n^{(0)}}{q_{n, n+1}}\sum_{k=i}^{n-1} {\widetilde F}_k^{(i)}.\end{equation}
Thus, to prove (\ref{e8.0}), it suffices to show that
$$\sum_{k=i}^{n-1} {\widetilde F}_k^{(i)}
=\prod_{i+1\le k\le n-1}\left[1+ \frac{{\tilde q}_k^{(0)}}{q_{k, k+1}}\right],\qquad n>i\ge 0.$$
This clearly holds when $n=i+1$. Suppose that it holds when
$n=\ell$, then
$$\aligned
\sum_{k=i}^{\ell} {\widetilde F}_k^{(i)}
&=\sum_{k=i}^{\ell-1} {\widetilde F}_k^{(i)}+ {\widetilde F}_{\ell}^{(i)}\\
&=\sum_{k=i}^{\ell-1} {\widetilde F}_k^{(i)}+\frac{{\tilde q}_{\ell}^{(0)}}{q_{\ell, \ell+1}}\sum_{k=i}^{\ell-1} {\widetilde F}_k^{(i)}\quad\text{\rm (by (\ref{e8.4}))}\\
&=\left[1+\frac{{\tilde q}_{\ell}^{(0)}}{q_{\ell, \ell+1}}\right]\sum_{k=i}^{\ell-1} {\widetilde F}_k^{(i)}\\
&=\prod_{i+1\le k \le \ell}\left[1+ \frac{{\tilde q}_{\ell}^{(0)}}{q_{k, k+1}}\right]\quad\text{\rm (by inductive assumption).}
\endaligned$$
Therefore, the required assertion holds for $n=\ell$ and it then holds for all $n>i$ by induction. We have thus proved the first assertion.

(b) By assumption, we have ${\tilde q}_n^{(k)}=q_{10}-c_1$ for every $k< n$. Hence, by $(\ref{e3.4-00})$ and $(\ref{e0.6-0})$, we obtain
$$\aligned
{\widetilde m}_n&=\frac{1}{q_{n, n+1}}\bigg(1+ {\tilde q}_1^{(0)}\sum_{k=0}^{n-1}{\widetilde m}_k\bigg),\qquad n\ge 1,\\
{\tilde d}_n&=\frac{1}{q_{n, n+1}}\bigg(1+ {\tilde q}_1^{(0)}\sum_{k=0}^{n-1}{\tilde d}_k\bigg),\qquad n\ge 1.
\endaligned
$$
As in the last proof, by using induction, we obtain the explicit expressions of $(\widetilde m_n)$ and $\big(\tilde d_n\big)$.

To study the divergence of $\sum_n m_n$, we adopt the

\noindent{\bf Kummer Test} Let $(u_n)$ and $(v_n)$ be two sequences of
positive numbers. Suppose that $\sum_0^{\infty} 1/v_n =\infty$ and
the limit $\kappa:=\lim_{n\to\infty} \kappa_n$ exists, where
$$\kappa_n =v_n \cdot \frac{u_n}{u_{n+1}}-v_{n+1}.$$
Then, the series $\sum u_n$ converges or diverges according to $\kappa>0$ or $\kappa<0$ respectively.

Set $v_n\equiv n$ and $u_n=m_n$:
$${m}_n=\frac{1}{q_{n, n+1}}
\prod_{0\le k \le n-1}\bigg[1+ \frac{q_{10}}{q_{k, k+1}}\bigg],\qquad n\ge 0.$$
Then
$$v_n\frac{u_n}{u_{n+1}}-v_{n+1}
= \frac{n(q_{n+1, n+2}-q_{n, n+1}-q_{10})}{q_{n, n+1}+q_{10}}-1.$$
Hence $\sum_n u_n<\infty $ if $\kappa'>1$ \big(resp.
$\sum_n u_n=\infty $ once $\kappa'<1$\big).
Clearly, $\sum_n m_n=\infty $ implies
$\sum_n {F}_n^{(0)}=\infty $.
Hence
$$d=\lim_{n\to\infty}\frac{d_n}{{F}_n^{(0)}}=\frac{1}{q_{01}}.$$
Furthermore,
$$\sup_{k\in E} \sum_{n=0}^k\big(F^{(0)}_n d-d_n\big)
=F^{(0)}_0 d=d<\infty.$$
This gives us the strong ergodicity by Proposition \ref{p5.2}.

We mention that Proposition \ref{pe} (with $0< c_i\equiv \lambda<q_{10}$) is also available for this
model.\hfill $\square$

\begin{rem}{\rm  For exponential ergodicity, the following
sufficient condition
\begin{equation}
M:=\sup_{n\ge 1} \bigg[\sum_{k=1}^{n-1}{F}_k^{(0)}\bigg]
\left[\sum_{j=n}^{\infty}\frac{1}{q_{j, j+1}{F}_j^{(0)}}\right]
<\infty,
\end{equation}
introduced in \cite{MZ}, is sufficient for Example \ref{ex8.1} but is not for Example \ref{ex8.2}.}
\end{rem}

\noindent{\it Proof}\quad It is obvious that $M<\infty$ iff
\begin{equation}
\varlimsup_{n\to\infty} \bigg[\sum_{k=1}^{n-1}{F}_k^{(0)}\bigg]
\left[\sum_{j=n}^{\infty}\frac{1}{q_{j, j+1}{F}_j^{(0)}}\right]
<\infty.
\end{equation}
For Example \ref{ex8.1}, because $q_{j, j+1}{F}_j^{(0)}$ is growing exponentially fast and so it is easy to check that $M<\infty$.
For Example \ref{ex8.2}, it suffices to consider $q_{n, n+1}=b (n+1)$ for some
$b>0$. By Kummer test, one may show that
$$\sum_{j=n}^{\infty}\frac{1}{q_{j, j+1}{F}_j^{(0)}}
=\infty$$
for suitable $b>0$ and then $M=\infty$. \hfill $\square$
\\[4mm]

\noindent\bf{\footnotesize Acknowledgements}\quad\rm
{\footnotesize The authors acknowledge the support by NNSFC (No. 11131003), SRFDP (No. 20100003110005),
the ``985'' project from the Ministry of Education in China,
the Fundamental Research Funds for the Central Universities,
and the Project Funded by the Priority Academic Program Development of Jiangsu Higher Education Institutions.}\\[4mm]

\noindent{\bbb{References}}

\noindent{\bbb Appendix. Key formulas used in the proofs}\\[0.1cm]

\noindent(A) Solution to the Poisson equation $\Omega g=Qg+cg$:
$$\aligned
g_n&=g_0+\sum_{0\le k\le n-1}\sum_{0\le j\leq k}\frac{{\widetilde F}_k^{(j)}(f_j-c_jg_0)}{q_{j,j+1}},\qquad n\ge 0.
\endaligned$$
(B) Three sequences.
\begin{itemize}
\item[(a)] $\widetilde F$-sequence:
\begin{gather}
\widetilde F_i^{(i)}=1,\qquad  \widetilde F_n^{(i)}=\frac 1{q_{n,n+1}}
\sum_{k=i}^{n-1}\tilde q_n^{(k)}\widetilde F_k^{(i)},\qquad n>i\ge 0,\nonumber\qquad\; \text{(\ref{ee0.01})}\end{gather}
where
$$\tilde q_n^{(k)}=q_n^{(k)}-c_n:=\sum_{j=0}^kq_{nj}-c_n,\qquad 0\le k<n.\qquad\qquad\qquad \text{(\ref{ee0.02})}$$
\item[(b)] $\widetilde m$-sequence:
$${\widetilde m}_0=\frac{1}{q_{01}},\quad  {\widetilde m}_n=\frac 1{q_{n,n+1}}\bigg(1+\sum_{k=0}^{n-1}{\tilde q}_n^{(k)}{\widetilde m}_k\bigg),\qquad n\ge 1.
\qquad\;\; \text{(\ref{e3.4-00})}$$
\item[(c)] $\tilde d$-sequence:
$$\tilde d_0=0,\qquad \tilde d_n=\frac 1{q_{n,n+1}}\Bigg(1+\sum_{k=0}^{n-1}\tilde q_n^{(k)}\tilde d_k\Bigg),
\qquad n\ge 1.\qquad\qquad \text{(\ref{e0.6-0})}$$
\end{itemize}
Representation of the three sequences:
$$
\aligned
&{\widetilde F}_i^{(i)}=1,\qquad  {\widetilde F}_n^{(i)}=\sum_{k=i+1}^n
\frac{{\widetilde F}_n^{(k)}{\tilde q}_k^{(i)}}{q_{k, k+1}},\quad n\ge i+1;\qquad \text{(\ref{e0.4-0})}\\
&{\tilde d}_n=\sum_{1\le k\le n}\frac{\widetilde F_n^{(k)}}{q_{k,k+1}},\quad\text{ (\ref{e0.6-1})}\qquad\qquad\quad
{\widetilde m}_n
=\sum_{k=0}^n\frac{{\widetilde F}_n^{(k)}}{q_{k,k+1}},\qquad n\ge 0.\quad\text{ (\ref{e3.4-0})}
\endaligned$$
Relation of the three sequences:
$${\widetilde m}_n= \frac{1}{q_{01}}{\widetilde F}_n^{(0)}+ {\tilde d}_n,\qquad n\ge 0.\qquad\qquad\qquad\qquad\qquad\qquad\qquad\qquad\;\text{ (\ref{e4.4-1})}$$
\end{document}